\pdfoutput=1

\documentclass{amsproc}
\usepackage{overpic}

\newtheorem{theorem}{Theorem}[section]
\newtheorem{lemma}[theorem]{Lemma}
\newtheorem{proposition}[theorem]{Proposition}
\newtheorem{corollary}[theorem]{Corollary}

\theoremstyle{definition}
\newtheorem{definition}[theorem]{Definition}

\newtheorem{conjecture}[theorem]{Conjecture}

\theoremstyle{remark}

\numberwithin{equation}{section}
\newcommand{\field}[1]{\mathbb{#1}}
\newcommand{\R}{\field{R}}
\newcommand{\Z}{\field{Z}}
\newcommand{\N}{\field{N}}
\newcommand{\C}{\field{C}}
\renewcommand{\P}{\field{P}}

\newcommand{\VV}{{\mathcal V}}

\newcommand{\RR}{{\mathcal R}}

\renewcommand{\AA}{{\mathcal A}}

\newcommand{\card}{\mathop{\rm card}}
\newcommand{\cp}{\mathrm{cap}}

\newcommand{\supp}{\mathop{\rm supp}}
\newcommand{\vt}[1]{\boldsymbol{#1}}

\newcommand{\const}{{\rm const}}

\renewcommand{\Re}{\mathop{\rm Re}}
\renewcommand{\Im}{\mathop{\rm Im}}
\newcommand{\dist}{\mathop{\rm dist}}
\newcommand{\isdef}{\stackrel{\text{\tiny def}}{=}}

\newcommand{\grad}{\mathop{\rm grad}}

\begin{document}

\title[Asymptotics of Heine-Stieltjes and Van Vleck polynomials]{On asymptotic behavior of Heine-Stieltjes and Van Vleck polynomials}

\author[A.~Mart\'{\i}nez-Finkelshtein]{Andrei Mart\'{\i}nez-Finkelshtein}
\address{Department of Statistics and Applied Mathematics
University of Almer\'{\i}a, SPAIN, and
Instituto Carlos I de F\'{\i}sica Te\'{o}rica y Computacional,
Granada University, SPAIN}
\email{andrei@ual.es}
\thanks{The first author was partially supported by Junta de Andaluc\'{\i}a, grants FQM-229, FQM-481, and P06-FQM-01735, as well as by the
research project MTM2008-06689-C02-01 from the Ministry of Science and Innovation of Spain and the European Regional Development Fund (ERDF)}

\author[E.A.~Rakhmanov]{Evgenii A.~Rakhmanov}
\address{Department of Mathematics,
University of South Florida, USA
}
\email{rakhmano@math.usf.edu}
\thanks{The second author was partially supported by  the NSF grant DMS-9801677.}

\date{February 3, 2009.}

\dedicatory{This paper is dedicated to the 60th Birthday of Guillermo (Bill) L\'opez Lagomasino}


\begin{abstract}
We investigate the strong asymptotics of Heine-Stieltjes polynomials -- polynomial solutions of a second order differential equations with complex polynomial coefficients. The solution is given in terms of critical measures (saddle points of the weighted logarithmic energy on the plane), that are tightly related to quadratic differentials with closed trajectories on the plane. The paper is a continuation of the research initiated in \cite{Martinez-Finkelshtein:2009ad}. However, the starting point here is the WKB method, which allows to obtain the strong asymptotics.
\end{abstract}

\maketitle


\section{Generalized Lam\'{e} equation}

In  1878 Heine studied the following problem: given two polynomials,
\begin{equation}\label{defAandB}
    A(z)=\prod_{i=0}^p (z-a_i)\,, \qquad B(z)=\alpha z^p+\text{lower
degree terms} \,, \quad \alpha \in \C ,
\end{equation}
where the zeros of  $A$ are complex and pairwise distinct, describe the polynomial solutions of the
\emph{generalized Lam\'{e} differential} equation (in algebraic form),
\begin{equation}
\label{DifEq}
A(z)\, y''(z)+ B(z)\, y'(z)- n(n+\alpha -1) V_n(z)\, y(z)=0,
\end{equation}
where $V_n$ is a polynomial (in general, depending on $n$) of degree $\leq p-1$; if $\deg V=p-1$, then $V$ is monic.
For $p=1$ we can easily recognize in \eqref{DifEq} the hypergeometric differential equation, while for $p=2$ it is known as the Heun's equation  (see \cite{Ronveaux95}).

Heine \cite{Heine1878} proved that for every $n \in \N$ there exist
at most
\begin{equation}\label{sigma}
\sigma(n)=\binom{n+p-1}{n}
\end{equation}
different polynomials $V_n$ such that \eqref{DifEq} admits a polynomial solution $y=Q_n \in \P_{n}$; hereafter $\P_n$ stands for the set of all algebraic polynomials of degree $\leq n$. These particular $V_n$'s are called \emph{Van Vleck polynomials}, and the corresponding
polynomial solutions $y=Q_n$ are known as \emph{Heine-Stieltjes} (or simply Stieltjes) polynomials. Furthermore, if the polynomials $A$ and $B$ are algebraically independent then for any $n\in \N$ there exist \emph{exactly} $\sigma(n)$ Van Vleck polynomials $V_n$, their degree is \emph{exactly} $p-1$, and for each $V_n$ equation \eqref{DifEq} has a unique (up to a constant factor) solution $y$ of degree $n$. According to recent results in  \cite{Shapiro2008a}, this statement still holds for arbitrary $A$ and $B$ for all sufficiently large $n$.

Stieltjes discovered an electrostatic interpretation of zeros of the polynomials discussed in \cite{Heine1878}, which attracted common attention to the problem.
He studied the problem \eqref{DifEq} in a particular setting, assuming that all zeros of $A$ are real and that all residues $\rho_k$ in \begin{equation}\label{BoverA}
    \quad \frac{B(x)}{A(x)}=
  \sum_{k=0}^p \frac{\rho_k}{x-a_k}
\end{equation}
are strictly positive.
He proved in \cite{Stieltjes1885} (see also \cite[Theorem
6.8]{szego:1975}) that in this case for each $n \in \N$ there are \emph{exactly}
$\sigma(n)$ different Van Vleck polynomials of degree $p-1$ and
the same number of corresponding Heine-Stieltjes polynomials $y$ of
degree $n$, given by all possible ways how the $n$ zeros of $y$
can be distributed in  the $p$ open intervals defined by the zeros of $A$.

After these two contributions a vast number of research papers on this topic has been published, dealing mainly with the real situation considered by Stieltjes. Besides the relevant contribution to the algebraic theory \cite{Shapiro2008a}, we can mention here the paper \cite{MR2003j:33031}, where the limit distribution of zeros of Heine-Stieltjes polynomials for the Stieltjes case has been established in terms of the traditional extremal problem for the weighted logarithmic energy on a compact set of the plane.

A new approach to the asymptotics of the polynomial solutions of \eqref{DifEq}, based on a different type of equilibrium problem, has been developed in \cite{Martinez-Finkelshtein:2009ad}, where zero distribution of Heine-Stieltjes polynomials was investigated. The present paper is in many senses a continuation of \cite{Martinez-Finkelshtein:2009ad}. We use essentially the same tools: families of continuous critical measures and their representations in terms of quadratic differentials with closed trajectories. However, the starting point of our considerations here is different: instead of electrostatics we use the WKB method, which allows to obtain the strong asymptotics.

\section{Asymptotic formula for Heine-Stieltjes polynomials}

In this section we describe briefly the strong asymptotics of polynomial solutions of \eqref{DifEq}. As we will see immediately, it is based on the asymptotics of Van Vleck polynomials, so formally we should have started with the latter topic. We choose another way of presentation, which follows rather the logic of the proof of the main result. We formulate next a theorem on asymptotics of Heine-Stieltjes polynomials, which contains Van Vleck  polynomials as unknown (``access'') parameters.

Hereafter we assume that $n\in \N$ is sufficiently large, so that there are exactly $\sigma(n)$ Van Vleck polynomials, each of degree exactly $p-1$, see \cite{Shapiro2008a}; denote by $V_{n}$ one of these Van Vleck  polynomials, and let $\VV_n$ stand for the set of zeros of $V_n$, $\AA=\{ a_0, \dots, a_{p} \}$ the set of zeros of $A$,  and  $\Omega_n \isdef  \C \setminus (\AA \cup \VV_n )$.

Fix any $z_0\in \C$; in any simply-connected domain $D\subset \Omega$ we can select a single-valued branch of $ \sqrt{V_{n}/A}$ and define the \emph{natural parameter}
\begin{equation}\label{defXi}
\xi=\xi_n(z)= \int_{z_0}^z  \sqrt{\frac{V_n(t)}{  A(t) }}\, dt\,.
\end{equation}
In what follows, many statements will be made in terms of trajectories of the quadratic differential
$$
\varpi_{n}=-\frac{V_{n}(z)}{A(z)}\, (dz)^2,
$$
which are basically level lines of the function $\Im \xi_{n}$.

We define also
\begin{equation}\label{zetaDefAlt}
\zeta_n(z)=\exp (\xi_n(z))= \exp \left( \int^z_{z_{0}}
\sqrt{\frac{V_n(t)}{ A(t) }}\, dt \right)\,,
\end{equation}
as well as the parameter
\begin{equation}\label{defLambda}
\lambda_n=n+\frac{\alpha-1}{2}\,.
\end{equation}

\begin{theorem}
\label{mainThm}
With the assumptions and notations above, let $y=Q_{n}$ be a Heine-Stieltjes polynomial (solution of \eqref{DifEq}) corresponding to $V_{n}$. Then:
\begin{enumerate}
	\item[(i)] there exist a set $\Gamma_{n} =\cup_{k\geq 1} \gamma_{k,n}$ comprised of a finite number of at most  $2p-1 $ disjoint Jordan arcs $\gamma_{k,n}$, such that for an appropriate selection of the branch in \eqref{defXi},
	\begin{equation}\label{asymptRepresentation3}
	y(z)=  H_n(z) \, \zeta_n^{\lambda_n}(z) \, \left(
	1+\varepsilon_n (z)\right)\,, \quad z \in \C\setminus \Gamma_{n},
	\end{equation}
	where
	\begin{equation}\label{defH}
	    H_n(z)=\left(\frac{A}{V_n}\, (z) \right)^{1/4} \, \exp
	    \left(-\int^z \frac{B}{2\, A}\, (t)\, dt
	    \right).
	\end{equation}
	Moreover, $\varepsilon_{n}(z)\to 0$ with $n\to \infty$ uniformly for $\dist(z, \Gamma_{n})\geq C>0$.

\item[(ii)] The endpoints of each $\gamma_{k,n}$ belong to $\AA \cup \VV_n$. Two arcs $\gamma_{j,n}$ and $\gamma_{k,n}$ for $i\neq j$ are either disjoint or have a common endpoint. The domain $\overline{\C}\setminus \Gamma_{n}$ is connected, and $\sqrt{V_{n}/A}$ has a single-valued branch. Consequently, the branch of the square root in \eqref{asymptRepresentation3} is fully determined by condition
$$
\lim_{z\to \infty} z \sqrt{\frac{V_n(z)}{ A(z) }}=1.
$$

\item[(iii)] Each $\gamma_{k,n}$ can be selected as a part of a trajectory of the quadratic differential
\begin{equation}\label{defQD}
\varpi_{n} = -\frac{V_{n}(z)}{A(z)}\, (dz)^2,
\end{equation}
which is close to its critical trajectory completed by two small segments connecting the trajectory with two points from $\AA \cup \VV_n$.

\end{enumerate}

\end{theorem}

Observe that formula \eqref{asymptRepresentation3} is given in terms of polynomials $V_{n}$. In this way, the problem of asymptotics of Heine-Stieltjes polynomials is reduced to that of Van Vleck  polynomials. 
We discuss the location of the zeros of $V_{n}$'s in Section \ref{sec:VanVleckzeros}.

The critical role  is played by the set $\Gamma_{n}=\cup_{k=1}^p \gamma_{k,n}$. More exactly, formula \eqref{asymptRepresentation3} should be combined with the local asymptotic formula for $Q_{n}$ in a neighborhood of $\gamma_{k,n}$ (but away from the endpoints):
\begin{equation}\label{Qn}
Q_n(z)=H_n(z) \, \left[  \zeta_n^{\lambda_n}(z) \, \left(
1+\varepsilon_{n,1}(z)\right)+   \zeta_n^{- \lambda_n}(z) \,
\left( 1+\varepsilon_{n,2}(z)\right)\right]\,,
\end{equation}
for certain small $\varepsilon_{n,j}(z)$. Local formulas of this type can be obtained directly from the WKB estimates as explained in \cite{Olver74b}. However, the existence of a global formula with the specified dominant term is not totally trivial. An arc $\gamma_{k,n}$ enters the WKB analysis naturally as the Stokes line: the dominant term in \eqref{Qn} changes across $\gamma_{k,n}$ (Stokes' phenomenon). This means that two terms are asymptotically equal on $\gamma_{k,n}$, and equation
$$
\zeta_n^{2\lambda_n}(z)=-1
$$
defines the zeros of $Q_{n}$.  In other words, arcs $\gamma_{k,n}$ asymptotically carry zeros of $Q_{n}$.

The proof of the theorem is based on the Liouvulle-Green (also known as the WKB) asymptotic formula. As we have mentioned above, the WKB method, the way it is explained in \cite{Olver74b}, is ready to use in order to obtain local estimates. Construction of a global formula from the local ones constitutes a challenge, as it usually happens in the WKB analysis. So, the proof of the theorem, which is outlined next, is not entirely on the surface, and relies heavily on the fact that both the coefficients and a solution of \eqref{DifEq} are polynomials.

\section{WKB analysis of the Heine-Stieltjes polynomials} \label{sec:WKB}

The Liouville-Green (also known as the WKB) method is a standard tool of asymptotic analysis of ordinary differential equations; we refer the reader to the comprehensive account contained in the monographs \cite{Fedoryuk:1993gf} and \cite{Olver74b}.

It was probably Nuttall \cite{MR891770} who applied the WKB analysis to the study of the strong asymptotics of the generalized Jacobi polinomials via a differential equation of the form \eqref{DifEq}.

\subsection{Liouville's transformation}

Recall that at this stage we assume that the Van Vleck  polynomial $V_n$ is fixed, and denote by $y=Q_n$ a corresponding Heine-Stieltjes polynomial of degree $n$. The following procedure is well known, so we summarize it here very briefly.

Rewriting equation \eqref{DifEq} in terms of
\begin{equation}\label{u}
u(z)=y(z) \exp \left( \int^z \frac{B}{2A}(t)\, dt\right)\,,
\end{equation}
we eliminate the first order derivative and get $ u''(z)=F_n(z) u(z)$, with
\begin{equation}\label{F}
F_n(z)=\left( \lambda_n^2 -\rho^2 \right)
\frac{V_n}{A}(z)+\left(\frac{B}{2A}
\right)^2(z)+\left(\frac{B}{2A} \right)'(z)\,,
\end{equation}
where $\lambda_{n}$ has been defined in \eqref{defLambda} and
$$
\rho=\frac{\alpha-1}{2}\,.
$$
Fix any $z_0\in \C$; the Liouville-Green transformation $w(\xi)=\left(\xi'\right)^{1/2} u$, with $\xi=\xi_{n}(z)$ defined in \eqref{defXi}, yields
\begin{equation}\label{DELiouville3}
     \frac{d^2 w}{d\xi^2} = \left(\lambda_n^2+ g_n(\xi )
    \right)\,w(\xi)\,,
\end{equation}
where
\begin{equation}\label{gStar}
g_n(\xi)=G_n(z(\xi))\,, \quad  G_n=\frac{A}{V_n}\, \left(
\left(\frac{B}{2\, A}\right)^2 +\left(\frac{B}{2\, A}\right)' -
T_n^2 -T_n' \right) -\rho ^2,
\end{equation}
and
$$
T_n(z)=\frac{1}{4}\,\left(\frac{A'}{A}-\frac{V'_n}{V_n}\right)(z)\,.
$$
A detailed discussion of this transformation in the complex plain and of the resulting WKB approximation can be found, e.g.\ in \cite{Olver74b}.

\subsection{Quadratic differentials on the Riemann sphere}

At this point let us recall briefly some notions from the theory of rational quadratic differentials on the Riemann sphere; see \cite{Pommerenke75} and \cite{Strebel:84}  for further details. Given the quadratic differential $\varpi_{n}$ defined in \eqref{defQD}, a smooth curve $\gamma$ along which
$$
\varpi_{n}=-V _{n}(z)/A(z)\, (dz)^2>0 \quad \Leftrightarrow \quad \Im \xi_{n}(z)=\const
$$
is a \emph{horizontal arc} of $\varpi_{n} $; the distinguished parameter $\xi_{n}$ was defined in \eqref{defXi}. More precisely, if $\gamma$ given by a parametrization $z(t)$, $t\in (\alpha ,\beta )$, then
$$
-\frac{V_{n}}{A}\, (z(t))\, \left(\frac{dz}{dt}\right)^2>0, \quad t
\in (\alpha ,\beta )\,.
$$
A maximal horizontal arc is called a \emph{horizontal trajectory} (or simply a \emph{trajectory}) of $\varpi_{n} $.
Analogously, trajectories of $-\varpi_{n}$ are called \emph{orthogonal} or  \emph{vertical trajectories} of  $\varpi_{n}$; along these curves
$$
V_{n}(z)/A(z)\, (dz)^2>0 \quad \Leftrightarrow \quad \Re \xi_{n}(z)=\const\,.
$$

We can define a conformal invariant metric associated with the quadratic differential $\varpi$, given by the length element $|d \xi_{n}|=|\sqrt{V_{n}/A}|(z) |dz|$; the $\varpi_{n}$-length of a curve $\gamma$ is
$$
\|\gamma\|_{\varpi_{n} }= \frac{1}{\pi}\, \int_{\gamma}
\sqrt{\left|\frac{V_{n}}{A}\right| } \,(z)\, |dz|\,;
$$
(observe that this definition differs by a normalization constant from the definition
5.3 in \cite{Strebel:84}). Furthermore, if $D$ is a simply connected domain not containing singular points of $\varpi_{n}$, we can introduce the $\varpi_{n}$-distance by
$$
\dist(z_1, z_2; \varpi, D)=\inf \{\|\gamma\|_{\varpi_{n}  }:
\, z_1, z_2 \in \bar \gamma, \; \gamma \subset D \}\,.
$$
Trajectories and orthogonal trajectories are in fact geodesics (in the $\varpi_{n}$-metric) connecting any two of its
points (see \cite[Thm. 8.4]{Pommerenke75}).

A simply connected domain $D$ not containing points from $\AA \cup \VV_{n}$  is called a \emph{$\varpi_{n}$-rectangle} if it is delimited by two horizontal and two vertical arcs of $\varpi_{n}$; in other words, if $\xi_{n}(D)$ is a (euclidean) rectangle $[a,b]\times [c,d]$, and $D\mapsto \xi_{n}(D)$ is a one-to-one conformal mapping. We call the value $d-c$ the $\varpi_n$-height,  $b-a$ the $\varpi_{n}$-length, and $2(b-a+d-c)$ the $\varpi_{n}$-perimeter of $D$. Obviously, these definitions are consistent with the freedom in the selection of the natural parameter $\xi_{n}$.

\subsection{Local asymptotics}

Theorems on local asymptotics of solutions are well-known (see \cite{Fedoryuk:1993gf}  and \cite{Olver74b}); using the local estimates that appear in \cite[Ch.\ VI, Theorem 11.1]{Olver74b} it follows that if $D$ is a $\varpi_n$-rectangle such that function $ g_{n}$ defined in (\ref{gStar}) is holomorphic in $\xi_n(D)=[a,b]\times [c,d]$, then
the differential equation (\ref{DELiouville3}) has in $[a,b]\times [c,d]$ two linearly independent holomorphic solutions $w_j$, $j=1,2$, of
the form
\begin{equation}\label{wwithxi11}
    w_j(\xi)= \exp \left\{ (-1)^{j+1} \lambda_n \xi \right\} \, \left(
1+\varepsilon _j(\xi)\right)\,,
\end{equation}
and such that
\begin{equation}\label{estimatesOlver}
|\varepsilon_j(\xi)|\leq \exp\left( \frac{1}{
|\lambda_n|} \, \int_{s_j}^\xi |g_n(t)|\, |dt| \right)-1\,, \quad \xi \in [a,b]\times [c,d]\,.
\end{equation}

The integrals here are taken following the \emph{progressive paths}, i.e.\ contours along which $\Re (\xi)$ is non-decreasing (for $j=1$) or non-increasing (for $j=2$). In the case of the rectangle we may take $s_1=a+ic$ and $s_2=b+id$, so that the whole rectangle is reachable by progressive paths. Taking advantage of the fact that the coefficients of the original equation \eqref{DifEq} are polynomials, we can estimate the total variation
$$
M(D)=\max_\gamma \int_\gamma \left| g_n(t) \right|\, |dt |\,,
$$
where $\gamma$ is any horizontal or vertical segment in $[a,b]\times [c,d]$. This yields the following result:
\begin{proposition} \label{prop:WKB1}
With the assumptions above, let $\widetilde{D}$ be an Euclidean rectangle such that $\xi_{n}(D)\subset \widetilde{D} $ and $\xi_{n}$ can be continued holomorphically to $\widetilde{D}$ in such a way that $\xi_{n}(\AA \cup \VV_{n})\cap \widetilde{D} = \emptyset$. Let  $d_n=\dist(\partial \widetilde D, \xi_{n}(D))$ be the Euclidean distance from $\xi_n(D)$ to the boundary of $\widetilde D$. Then
 we can replace the estimates \eqref{estimatesOlver} by
\begin{equation}\label{boundsOlver}
|\varepsilon_j(\xi)|\leq \exp\left( \frac{K}{ d_n\,
|\lambda_n|} \right)-1\,, \quad \xi \in \RR\,.
\end{equation}
\end{proposition}

In a $\varpi_{n}$-rectangle $D$ we can select a single valued branch of the function $H_{n}$ introduced in \eqref{defH}.
Then a direct consequence of the proposition above is
\begin{corollary} \label{cor:WKB1}
Let $D$ be a $\varpi$-rectangle. Then a general solution of
(\ref{DifEq}) in $D$ has the form
\begin{equation}\label{generalSol1}
y(z)= H_n(z) \, \left[\kappa_1\, \zeta_n^{\lambda_n}(z) \, \left(
1+\varepsilon _1(z)\right)+ \kappa_2\, \zeta_n^{- \lambda_n}(z) \,
\left( 1+\varepsilon _2(z)\right)\right]\,,
\end{equation}
with $\zeta_{n}$ defined in \eqref{zetaDefAlt}.
We have
\begin{equation}\label{boundOn G}
|\varepsilon_j(z)|\leq \exp\left( \frac{K}{ d_n \,
|\lambda_n|} \right)-1\,,  \quad z \in
D\,,
\end{equation}
where $d_n$ is the Euclidean distance defined in Proposition \ref{prop:WKB1}.
\end{corollary}

\subsection{Global asymptotic formula away from zeros}

The result above shows that if we stay away from the singularities $\AA \cup \VV_{n}$, we can control the errors in the WKB approximation uniformly. This motivates the following definition.
For $\varepsilon >0$ and $t \in \C$ and subset $K \subset \C$ we
denote
$$
D_\varepsilon(t )\isdef \{z\in \C:\, |z-t|<\varepsilon  \}\,, \quad
D_\varepsilon(K ) \isdef \bigcup_{t\in K} D(t, \varepsilon )\,,
$$
and
$$
D_{n,\varepsilon } \isdef D_\varepsilon(\AA  \cup  \VV_n)\,.
$$

Let $D$ be a $\varpi_n$-rectangle in $\C \setminus D_{n,\varepsilon }$ containing infinity. We can select the branch of $\zeta_n$ and the initial point $z_{0}$ in \eqref{zetaDefAlt} in such a way that $|\zeta_n(z)|\geq 1+\delta$ for a $\delta>0$, and there is a constant $\kappa$ such that
\begin{equation}\label{QnOld}
Q_n(z)=H_n(z) \, \left[  \zeta_n^{\lambda_n}(z) \, \left(
1+\varepsilon _1(z)\right)+   \kappa\, \zeta_n^{- \lambda_n}(z) \,
\left( 1+\varepsilon _2(z)\right)\right]\,.
\end{equation}
It follows that if we continue all functions analytically in the maximal domain containing  $D$ such that $|\zeta_n(z)|\geq 1+\delta$, then formula \eqref{QnOld} still holds.

Next, standard arguments show that in a domain $D$, if $Q_n\neq 0$ for all
sufficiently large $n$, then there exists a dominant term in
(\ref{generalSol1}):
\begin{lemma} \label{lemma:nozeros}
Let $D=\{ z\in \C:\, a <\Re (\xi_n)<b, \; c <\Im (\xi_n)<
d\}$ be a simply connected $\varpi _n$-rectangle in $\C \setminus
D_{n,\varepsilon }$, with
$$
d-c > \frac{\pi}{|\lambda_n|}\,.
$$
If $\Omega$ does not contain zeros of a solution $y$ of
(\ref{DifEq}) then for any choice of the lower limit of
integration and with an appropriate choice of the branch of the
square root in (\ref{zetaDefAlt}) we have
\begin{equation}\label{asymptRepresentation1}
y(z)= H_n(z) \, \zeta_n^{\lambda_n}(z) \, \left(
1+\varepsilon_n (z)\right)\,.
\end{equation}
Moreover, there exists a constant $M=M(\varepsilon , D)$,
independent of $n$, such that  for any $\varpi_n$-rectangle
$D'\subset D$,
\begin{equation}\label{boundM}
|\varepsilon_n(z)|\leq \frac{M}{ |\lambda_n|} \,, \quad t \in
D'\,.
\end{equation}
\end{lemma}

Now let us observe what happens if we have a union of two
$\varpi_n$-rectangles, $D^{(1)}$ and $D^{(2)}$, such
that $D^{(1)} \cap D^{(2)}$ has an interior point,
$z_0$. Assume first that we take this $z_0$ as the lower limit of
integration in (\ref{zetaDefAlt}) and (\ref{defH}). From Lemma
\ref{lemma:nozeros} it follows that if $y \not \equiv 0$,
$$
y(z)=\kappa^{(j)} H_n(z) \, \zeta_n^{\lambda_n}(z) \, \left(
1+\varepsilon_n^{(j)} (z)\right)\,, \quad z\in D^{(j)}\,,
\quad j = 1, 2,
$$
for certain non zero constants $\kappa^{(1)}$ and $\kappa^{(2)}$.
Evaluating at $z=z_0$ we get
$$
\kappa^{(1)} \, \left( 1+\varepsilon_n^{(1)}
(z_0)\right)=\kappa^{(2)} \, \left( 1+\varepsilon_n^{(2)}
(z_0)\right)\,, \quad z\in D^{(1)} \cap D^{(2)}\,,
$$
so that,
$$
\kappa^{(2)}=\kappa^{(1)} \, \frac{ 1+\varepsilon_n^{(1)} (z_0)}{
1+\varepsilon_n^{(2)} (z_0)}=\kappa^{(1)} \, \left(
1+\varepsilon_n (z)\right)\,, \quad |\varepsilon_n(z)|\leq
\frac{M'}{ |\lambda_n|} \,.
$$
In the intersection (neighborhood of $z_0$) both expressions
for $y$ should match, hence we have that in $D^{(1)} \cap
D^{(2)}$,
$$
\zeta_n^{\lambda_n}(z) \, \left( 1+\varepsilon_n^{(1)}
(z)\right)=\zeta_n^{\pm \lambda_n}(z) \, \left( 1+\varepsilon_n
(z)\right)  \,  \left( 1+\varepsilon_n^{(2)} (z)\right)\,.
$$
Since $D^{(1)} \cap D^{(2)}$ is an open set, necessarily
the same branch of $\zeta_n$ has been taken in both sides of the
previous identity. This argument shows that Lemma \ref{lemma:nozeros} is valid in
the union $D^{(1)} \cup D^{(2)}$, eventually with a
different constant in the right hand side of (\ref{boundM}).

Now, modifying the lower limit of integration only changes the
normalization constant $\kappa$ in (\ref{asymptRepresentation1}).
That means that in the representation
(\ref{asymptRepresentation1}) we can choose any lower bound $z_0
\in D^{(1)} \cup D^{(2)}$, as long as we understand the
function $\zeta_n$ as the analytic continuation along a path in
$D^{(1)} \cup D^{(2)}$.

\begin{figure}[htb]
\centering \mbox{\begin{overpic}[scale=0.65]%
{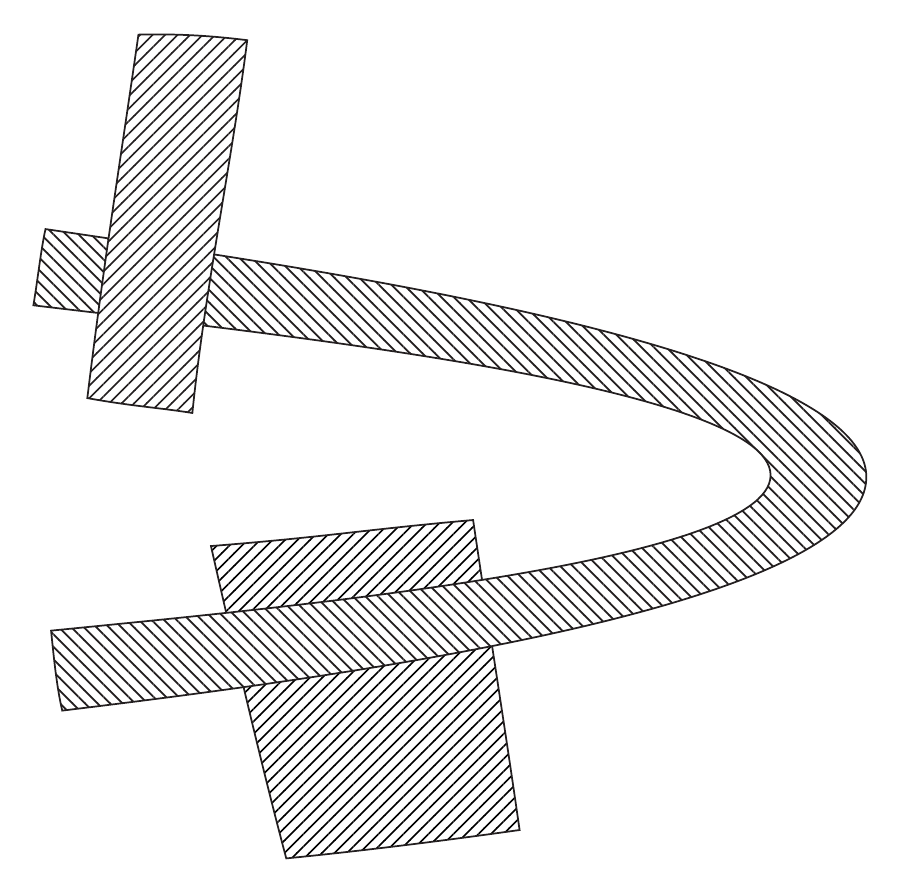}%
\put(15,85){$D^{(3)} $}
\put(40,14){$D^{(1)} $}
\put(8,23){$D^{(2)} $}
\end{overpic}}
\caption{A $\varpi_{n}$-chain. } \label{fig_chain}
\end{figure}

The discussion above motivates the following definition:
\begin{definition}
Let $D^{(j)}$, $j=1, \dots, k$, be a finite set of
$\varpi_n$-rectangles with $\varpi_n$-height greater than
$|\lambda_n|$ and bounded $\varpi_n$-perimeters. If each
intersection $D^{(j)} \cap D^{(j+1)}$, $j=1, \dots,
k-1$, has an interior point, then their union $\cup_{j=1}^k
D^{(j)}$ is called a finite \emph{$\varpi_n$-chain} (see Figure \ref{fig_chain}).
\end{definition}

Thus,  we have proved the following result:
\begin{lemma} \label{lemma:nonzerosGeneral}
Let $\widetilde D\subset \C \setminus D_{n,\varepsilon }$ be a domain, and
$D =\cup_{j=1}^k D ^{(j)}\subset \widetilde D$ a finite and simply connected
$\varpi_n$-chain in $\widetilde D$. If a solution $y$ of (\ref{DifEq})
does not vanish in $D$ then for any choice of the lower limit
of integration in $D$, and with an appropriate choice of the
branch of the square root in (\ref{zetaDefAlt}) there exists a
constant $C$ such that
(\ref{asymptRepresentation1})--(\ref{boundM}) holds in $D$.
The branch of $\zeta_n$ is obtained by analytic continuation in
$D$.
\end{lemma}

In other words, asymptotic representation can be continued along
the chains of $\varpi_n$-rectangles, as long as we stay away from
the zeros of the solution and of the sets $\AA$ and $\VV_n$.

\subsection{Zeros of a solution of the differential equation}

Assume now that $y$ is a nontrivial solution of (\ref{DifEq})
and $z_0\in \C \setminus D_{n,\varepsilon }$ is a zero of $y$. Let
$\Omega$ be a maximal simply connected $\varpi_n$-rectangle in
$\C \setminus D_{n,\varepsilon }$ containing $z_0$, and let
$\gamma$ be the vertical trajectory in $\gamma$ passing through
$z_0$. By assumption, $z_0$ is a regular point of the quadratic
differential $\varpi_n$, so that $\gamma$ is well defined.

In $\Omega$ the expression (\ref{generalSol1}) is valid. In
particular,
\begin{equation}\label{characterizationZero0}
y(z)=0 \quad \Leftrightarrow \quad 2\lambda_n\, \xi_n(z)= \log
\left(-\frac{\kappa_2}{\kappa_1}\,
\frac{1+\varepsilon_1(z)}{1+\varepsilon_2(z)}\right) \quad \mod
(2\pi i)\,.
\end{equation}
Hence, from the assumption $y(z_0)=0$ it follows that
$$
2\lambda_n\, \xi_n(z_0)= \log \left(-\frac{\kappa_2}{\kappa_1}\,
\frac{1+\varepsilon_1(z_0)}{1+\varepsilon_2(z_0)}\right) \quad
\mod (2\pi i)\,.
$$

Assume that $z_1$ is another point on $\gamma$ that satisfies the
following condition:
\begin{equation}\label{oneStep}
    \frac{\lambda_n}{\pi}\, \int_{z_0}^{z_1}\left| d\xi_n(t)
    \right|= \lambda_n \left\| \gamma(z_0, z_1)
    \right\|_{\varpi_n}\in \N\,,
\end{equation}
where $\gamma(z_0, z_1)$ is the arc of $\gamma$ joining $z_0$ and
$z_1$.

\begin{definition}
Let $z_0\in \C \setminus D_{n,\varepsilon }$, and let $\gamma$ be
the largest connected vertical arc in $\C \setminus
D_{n,\varepsilon }$ of the quadratic differential $\varpi_n$
passing thorough $z_0$. The set
$$
\omega=\bigcup_j \omega _j
$$
is a \emph{necklace} in $\C \setminus D_{n,\varepsilon }$
corresponding to $z_0$ if all $\omega_j$, called \emph{beads}, are
$\varpi_n$-rectangles of the form
$$
\omega_j=\{z\in \Omega:\, |\Re (\xi_n(z)-\xi_n(z_j))|<\delta, \,
|\Im (\xi_n(z)-\xi_n(z_j))|<\delta\}\,,
$$
where each $z_j$ satisfies condition (\ref{oneStep}), there exists a constant $M$ such that $\lambda_{n}^2 \delta\leq M$, and all $\omega_j \subset \C \setminus
D_{n,\varepsilon }$. The vertical arc $\gamma$ is the
\emph{string} of the necklace.
\end{definition}
\begin{figure}[htb]
\centering \includegraphics[scale=0.8]{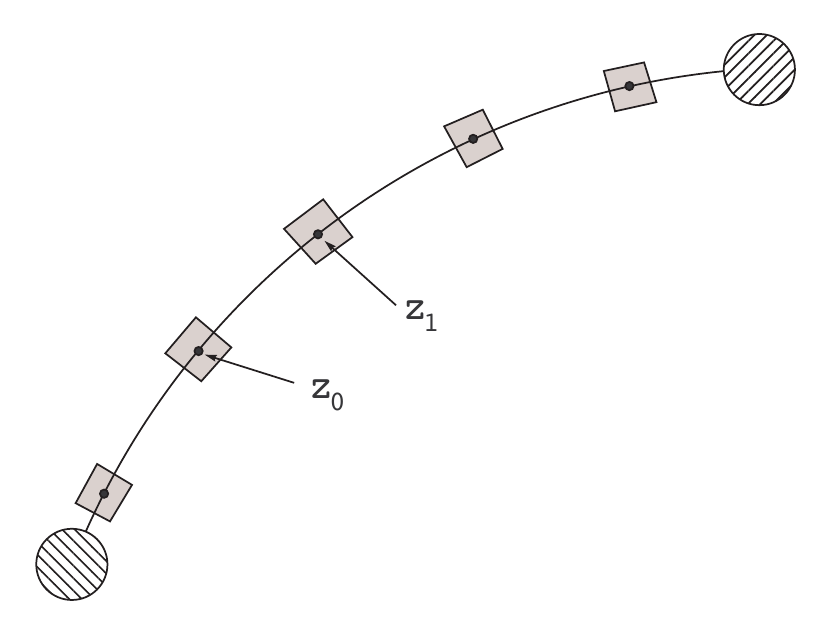} \caption{Necklace
corresponding to $z_0$. } \label{fig_necklace}
\end{figure}

A direct consequence of Rouche's theorem is the following statement:
\begin{lemma} \label{necklace}
Let $z_0\in \C \setminus D_{n,\varepsilon }$ be a zero of a
nontrivial solution $y$ of (\ref{DifEq}), and let $\omega$ be a
necklace in $\C \setminus D_{n,\varepsilon }$ corresponding to
$z_0$. Then each bead $w_j\subset \C \setminus D_{n,\varepsilon }$
contains one and only one zero of $y$.

Furthermore, let  $\Omega$ be a simply connected $\varpi _n$-rectangle in $\C
\setminus D_{n,\varepsilon }$, and $z_0\in \Omega$ be a zero of a
nontrivial solution $y$ of (\ref{DifEq}). If  $\omega$ is the
necklace in $\C \setminus D_{n,\varepsilon }$ corresponding to
$z_0$, then $y(z)\neq 0$ for $z \in \Omega \setminus \omega$.

\end{lemma}

These results   are applicable to any solution of (\ref{DifEq}).
Now we concentrate on the $n$-th degree Heine-Stieltjes polynomial. The main difference is that we know that it has exactly $n$ zeros, with account of multiplicity. This and Lemma \ref{necklace} immediately yield the following
\begin{proposition}
\label{prop:neclaceHS} Assume that $z_0\in \C \setminus
D_{n,\varepsilon }$ is a zero of a Heine-Stieltjes polynomial
$Q_n$, and let $\gamma$ be the largest connected vertical arc in
$\C \setminus D_{n,\varepsilon }$ of the quadratic differential
$\varpi_n$ passing thorough $z_0$. Then
$$
\|\gamma\|_{\varpi_n}\leq \frac{n+2}{\lambda_n}\,.
$$
In particular, every string of a zero-carrying necklace has a
finite $\varpi_n$-length, and starts and ends at $D_{n,
\varepsilon }$, or is a closed curve.
\end{proposition}
Indeed, by Lemma \ref{necklace}, every necklace carries at least
$\lambda_n \|\gamma\|_{\varphi_n}-2$ zeros of $Q_n$, where
$\gamma$ is the string of the necklace.

Summarizing, if $Q_n$ is a Heine-Stieltjes polynomial, there exist
a finite number of zero-carrying necklaces with the corresponding
strings $\gamma_{n,j}$ such that all the zeros belong to the union
of these necklaces or lie in $D_{n, \varepsilon }$.

\section{Asymptotics of Van Vleck polynomials} \label{sec:VanVleckzeros}

As it is clear from Theorem \ref{mainThm}, the zeros of the Van Vleck polynomial $V_{n}$ are the key parameters in the asymptotic expression for the Heine-Stieltjes polynomials. Hence, as a next step we derive a set of equation that will characterize their positions. We start with a formal argument and postpone to next section a more detailed discussion about consistency and meaning of these equations.

\begin{theorem}
\label{thm:VVzeros}
With the assumptions and notations above, let $y=Q_{n}$ be a Heine-Stieltjes polynomial (solution of \eqref{DifEq}) corresponding to $V_{n}$, and let $\gamma_{n,k}$ be the set of arcs defined in Theorem \ref{mainThm}. If there exists an $ \varepsilon>0 $ such that all arcs $\gamma_{k,n}$ are disjoint and  the $\varpi_{n}$-distance between them is $> \varepsilon$, then the following system of equations is satisfied:
\begin{equation}\label{VVequationsGeneral}
\frac{1}{\pi i}\, \int_{\gamma_{k,n}} \sqrt{\frac{V_{n}(t)}{A(t)}}\, dt = \frac{  m_{k} }{\lambda_{n}}+ \frac{ \rho_{k}/2 +   \delta_{k,n}}{\lambda_{n}}, \quad k=1, \dots, p-1, \quad m_{k}\in \mathbb M_{n},
\end{equation}
where the index set $\mathbb M_{n}$ is a finite subset of $\N\cup \{ 0\}$.
If for $\gamma_{k,n}$ we denote by $\eta_{k}$  the set of endpoints of $\gamma_{k,n}$ that belong to $\AA$, then
$$
\rho_{k}\isdef  \sum_{a\in\eta_{k}} \frac{B(a)}{A'(a)} + 1 - \card (\eta_{k}).
$$
Furthermore, there exists a constant $C=C(\varepsilon)$ such that
$$
|\delta_{k,n}|\leq \frac{C}{n}.
$$
 \end{theorem}

This theorem is a simple consequence of the asymptotic formula \eqref{asymptRepresentation3}, which is valid in a neighborhood of $\gamma_{k,n}$. Observing its increment  along a closed Jordan curve $\widetilde \gamma_{k,n}$  encircling $\gamma_{k,n}$ in the positive direction and using the argument principle we get the formulas above.

The system of equation \eqref{VVequationsGeneral} consists of $p-1$ equations (since we have fixed the residue at infinity, equation $k=p$ is dependent from the other $p-1$ ones). More exactly,  \eqref{VVequationsGeneral} presents a collection of systems of equations. In order to define it completely we need to specify:
\begin{enumerate}
\item[(a)] the combinatorics: points from $\AA \cup \VV_{n}$ are arranged in $p$ pairs (endpoints of $\gamma_{k,n}$).
\item[(b)] the homotopic types of curves $\gamma_{k,n}$ (once the combinatorics is fixed); and
\item[(c)] the range of values the integer parameters $m_{k}$ in the right hand side of  \eqref{VVequationsGeneral} may take.

\end{enumerate}

Observe that by  \eqref{VVequationsGeneral},
$$
\frac{1}{\pi i}\, \int_{\gamma_{k,n}} \sqrt{\frac{V_{n}(t)}{A(t)}}\, dt = \frac{  m_{k} }{\lambda_{n}}+ \frac{ \rho_{k}/2  }{\lambda_{n}} +\mathcal O\left( \frac{1}{n^2}\right), \quad k=1, \dots, p-1, \quad m_{k}\in \N.
$$
This system allows eventually  to to find the position of ``almost all'' zeros of Van Vleck polynomials, with an error smaller than the distance between zeros. In most part of the range this error is $\mathcal O(n^{-2})$. A complete analysis however is cumbersome and contains a combination of analytic, geometric and combinatorial arguments. We restrict our presentation here to the case $p=2$ (three points), which at least has a trivial combinatorics and a rather simple geometry, with additional remarks  on the case $p=3$.

\subsection{Case of $p=2$} \label{subs:p=2Cheb}

We have $\AA=\{ a_{0}, a_{1}, a_{2}\}$ (in general, non-collinear) and want to discuss the issues (a)-(c) raised above in order to define completely the set of equations determining the position of the zero of the Van Vleck polynomial. The first ingredient we need is the solution of the classical minimal capacity problem posed in the class of all continua in $\C$ containing $\AA$ (the Chebotarev's problem). It was proved by Gr\"otzsch \cite{Grotzsch1930a} and Lavrentiev \cite{Lavrentieff1930, Lavrentieff1934} that there exists a unique $\Gamma^*=\Gamma^*(\AA)$ satisfying
\begin{equation}\label{defChebotarev}
 \cp(\Gamma^*) =\min \{ \cp(F):\, F \text{ a continuum containing } \AA\},
\end{equation}
where $\cp(\cdot)$ denotes the logarithmic capacity.
If $a_{j}$'s are not collinear, then  $\Gamma^*$ is a union of three arcs, $\gamma_{0}^{*}$, $\gamma_{1}^{*}$, and $\gamma_{2}^{*}$,  connecting a point $v^*=v^*(\AA)$ with points $a_{j}$'s, respectively (see also \cite{Kuzmina1980} and \cite{Ortega-Cerda2008}). We call this $\Gamma ^*$ the \emph{Chebotarev's compact} or \emph{Chebotarev's continuum} corresponding to $\AA$, and point $v^{*}$ is the \emph{Chebotarev's center} of the set $\AA$. If we define
\begin{equation}\label{defMuK}
M_{j}=\frac{1}{\pi}\, \int_{\gamma_{k}}\left| \frac{t-v^*}{A(t)}\right|^{1/2}|dt|, \quad k=0,1,2,
\end{equation}
then
$$
M_{0}+M_{1}+M_{2}=1.
$$
Observe that each $M_{j}$ is the $\varpi^*$-length of $\gamma_{j}^{*}$, where
$$
\varpi^*\isdef \frac{v^*-z}{A(z)}\, (dz)^2.
$$

Next we define three analytic functions (elements) $w_{k}(v)$, $k=0, 1, 2$. We describe first $w_{0}$ as a germ of an analytic function at $v=a_{0}$, which allows unlimited analytic continuation to $\C\setminus \AA$. For all $v$ in a sufficiently small neighborhood of $a_{0}$ let $\Delta_{0}$ be the segment $[a_{0}, v]$, and $\Delta_{1}=[a_{1},a_{2}]$. Denote $\Omega\isdef \C\setminus (\Delta_{0}\cup \Delta_{1})$. With $R(z)\isdef (z-v^*)/A(z)$ we consider in $\Omega$ the single-valued branch of $\sqrt{R}$ given by the asymptotic condition $\lim_{z\to \infty}z\sqrt{R(z)}=1$, and define
$$
w_{0}(v)=\frac{1}{\pi i}\, \int_{\Delta_{0}} \sqrt{R(t)}\, dt=\frac{1}{2 \pi i}\, \oint_{\partial \Delta_{0}} \sqrt{R(t)}\, dt,
$$
where $\partial \Delta_{0}$ is the doubly-connected component of $\partial \Omega$ (the Carath\'{e}odory boundary of $\Omega$), with the boundary values of the branch of $\sqrt{R}$ specified above. It is clear that $w_{0}$ is analytic at $v=a_{0}$, and may be continued as a multi-valued analytic function to $\C\setminus \AA$. Similarly, we define $w_{1}$ and $w_{2}$, starting from $v=a_{1}$ and $v=a_{2}$, respectively. The proof of the following statement can be found, e.g.~in \cite{Martinez-Finkelshtein:2009ad}:
\begin{proposition}
\label{prop:VVlocations}
For $k=0, 1, 2$ there exist an analytic arc $\ell_{k}$, connecting $a_{k}$ with $v^*$, defined by
$$
\ell_{k}=\{ v\in \C:\, w_{k}(v)\in [0, M_{k}]\}.
$$
Furthermore, for $k=0,1,2$, $w_{k}(v)$ is univalent in a neighborhood $\mathcal L_{k}$ of $\ell_{k}$ (which is mapped by $w_{k}$ onto a neighborhood of the interval $[0,M_{k}]$).
\end{proposition}

Now we are ready to formulate a theorem about the asymptotic location of the zeros of Van Vleck polynomials for $p=2$.
\begin{theorem}
\label{thm:VVzerosP=2}
Let $y=Q_{n}$ be a Heine-Stieltjes polynomial (solution of \eqref{DifEq}) with $\AA=\{ a_{0}, a_{1}, a_{2}\}$,  corresponding to $V_{n}(z)=z-v_{n}$. Assume that there exists an $ \varepsilon>0 $ such that $\dist(v_{n},v^*)>\varepsilon$ (where $v^{*}$ is the Chebotarev's center). Then for all sufficiently large $n$, there exist an  index $k\in \{0, 1, 2 \}$ and an integer value $0\leq m_{k} \leq [n  M_{n}] $ such that the following equation is satisfied:
\begin{equation}\label{VVequationP=2}
w_{k}(v_{n}) = \frac{  m_{k} }{\lambda_{n}}+\frac{1}{\lambda_{n}}\,  \frac{B(a_{k})}{2 A'(a_{k})}+  \frac{     \delta_{k,n}}{\lambda_{n}},
\end{equation}
and there exists a constant $C=C(\varepsilon)$ such that
$$
|\delta_{k,n}|\leq \frac{C}{n}.
$$
 \end{theorem}

Obviously, this statement does not cover, roughly speaking, $\varepsilon n$ zeros out of $n+1$ possible zeros of Van Vleck  polynomials, but it provides an error estimate of order $C(\varepsilon)n^{-2}$. The exceptional set of ``missed'' Van Vleck zeros may be made smaller (up to a constant) for the price of relaxing the error estimate; also the technical details become more cumbersome.

\subsection{Case of $p=3$}

Now we turn to the problems (a) and (b) related to the system \eqref{VVequationsGeneral}, namely, we discuss the combinatorics and the homotopic type of curves $\gamma_{k,n}$ from Theorem \ref{thm:VVzeros}; we leave the issue (c) of the range of values of $m_{k}$ in the equations to the following Section.

\begin{figure}[htb]
\centering \begin{tabular}{l@{\hspace{-1mm}}c} 
\hspace{-1.4cm}\mbox{\begin{overpic}[scale=0.4]%
{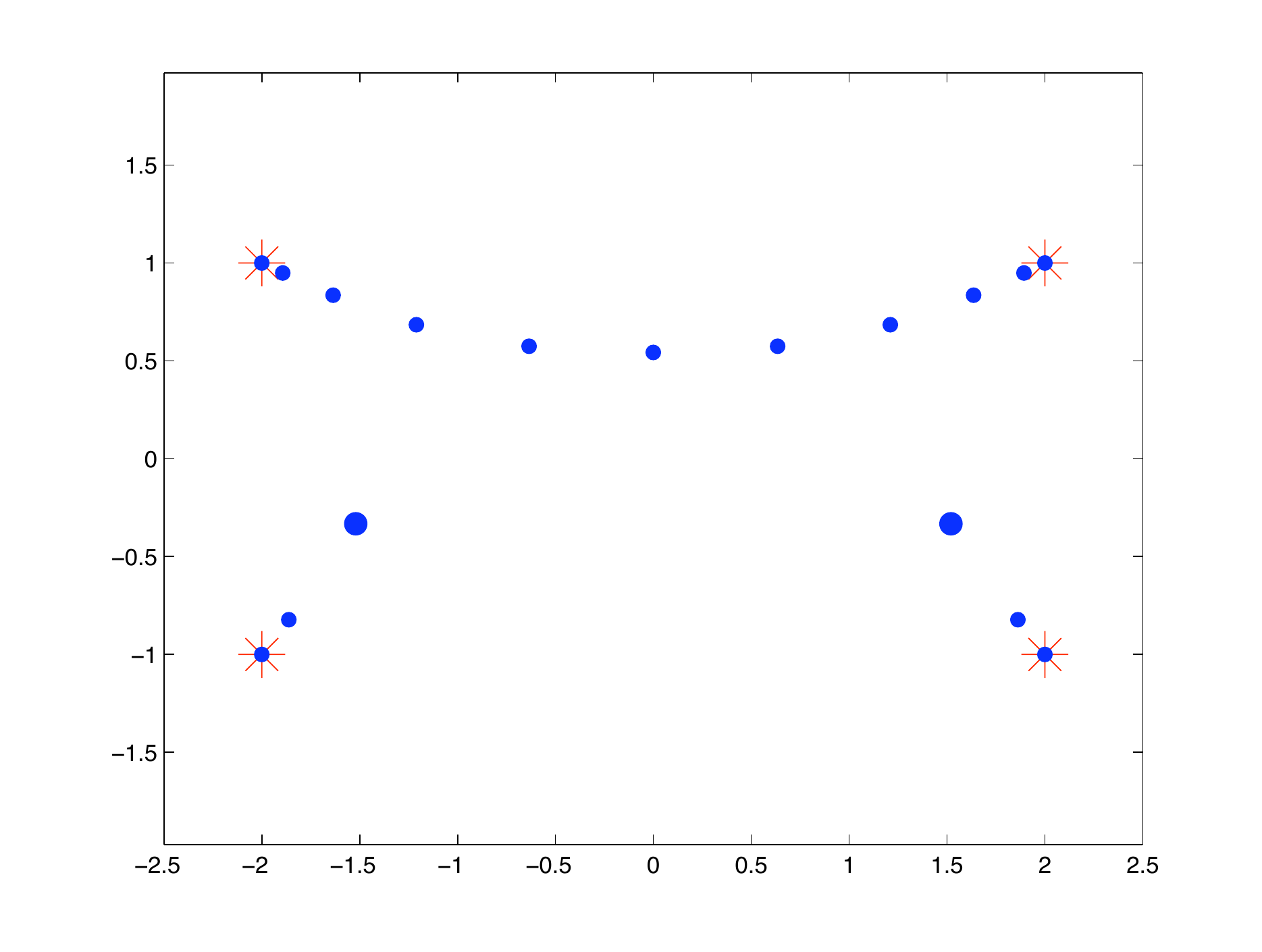}%
\end{overpic}} &
\hspace{-0.5cm}
\mbox{\begin{overpic}[scale=0.4]%
{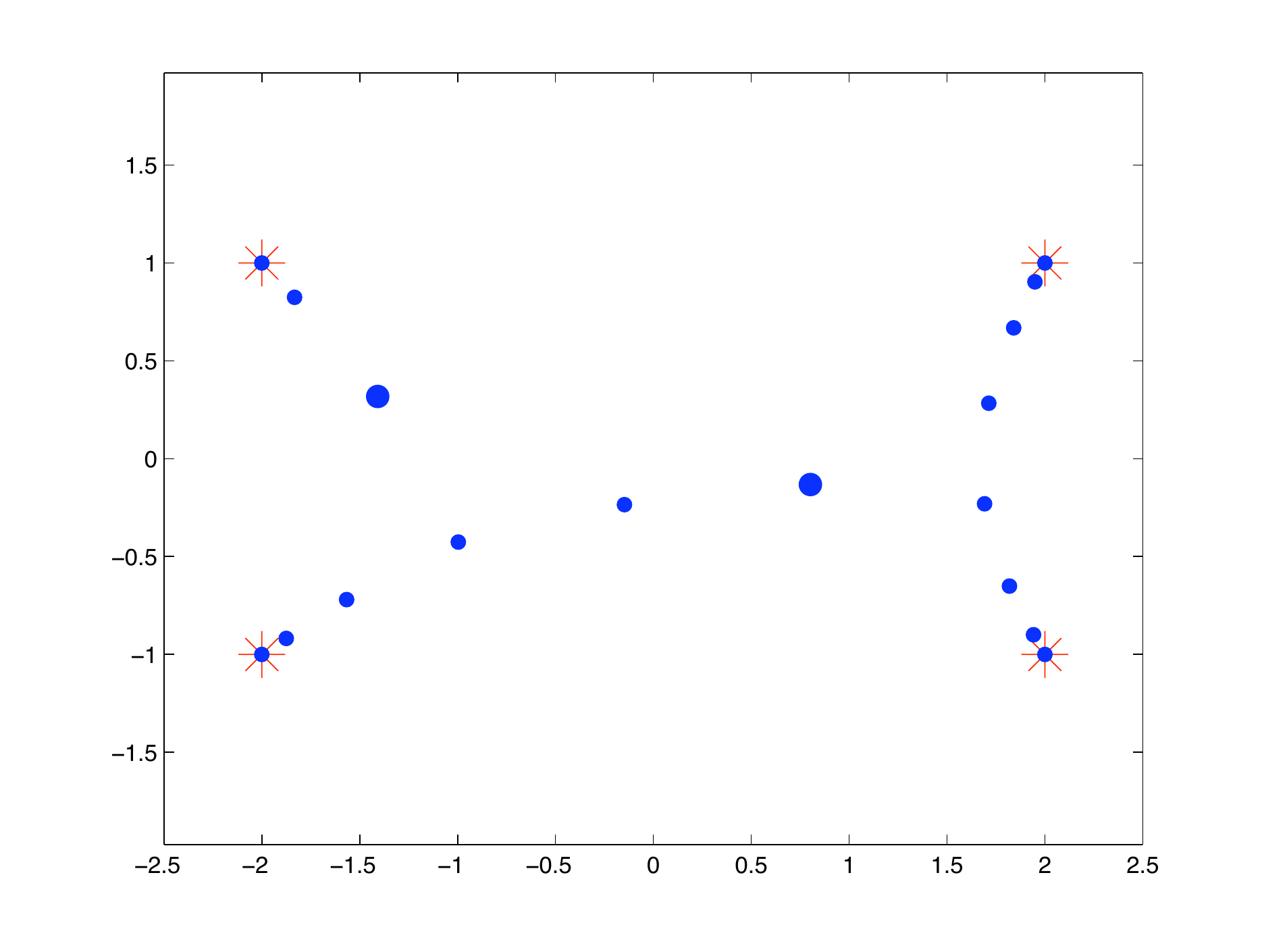}%
\end{overpic}}
\\[-1mm] 
\hspace{-1.4cm}\mbox{\begin{overpic}[scale=0.4]%
{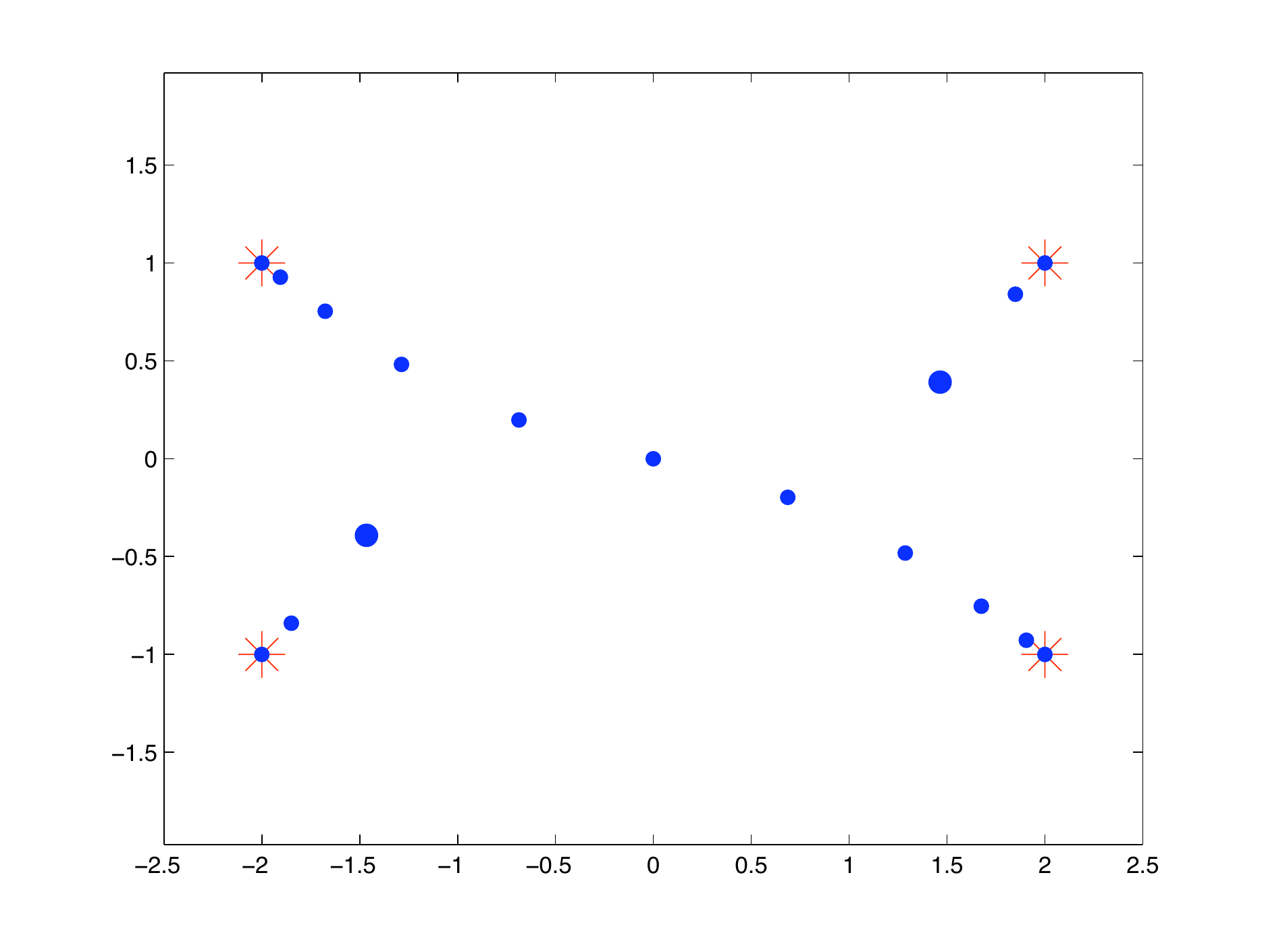}%
\end{overpic}} &
\hspace{-0.5cm}
\mbox{\begin{overpic}[scale=0.4]%
{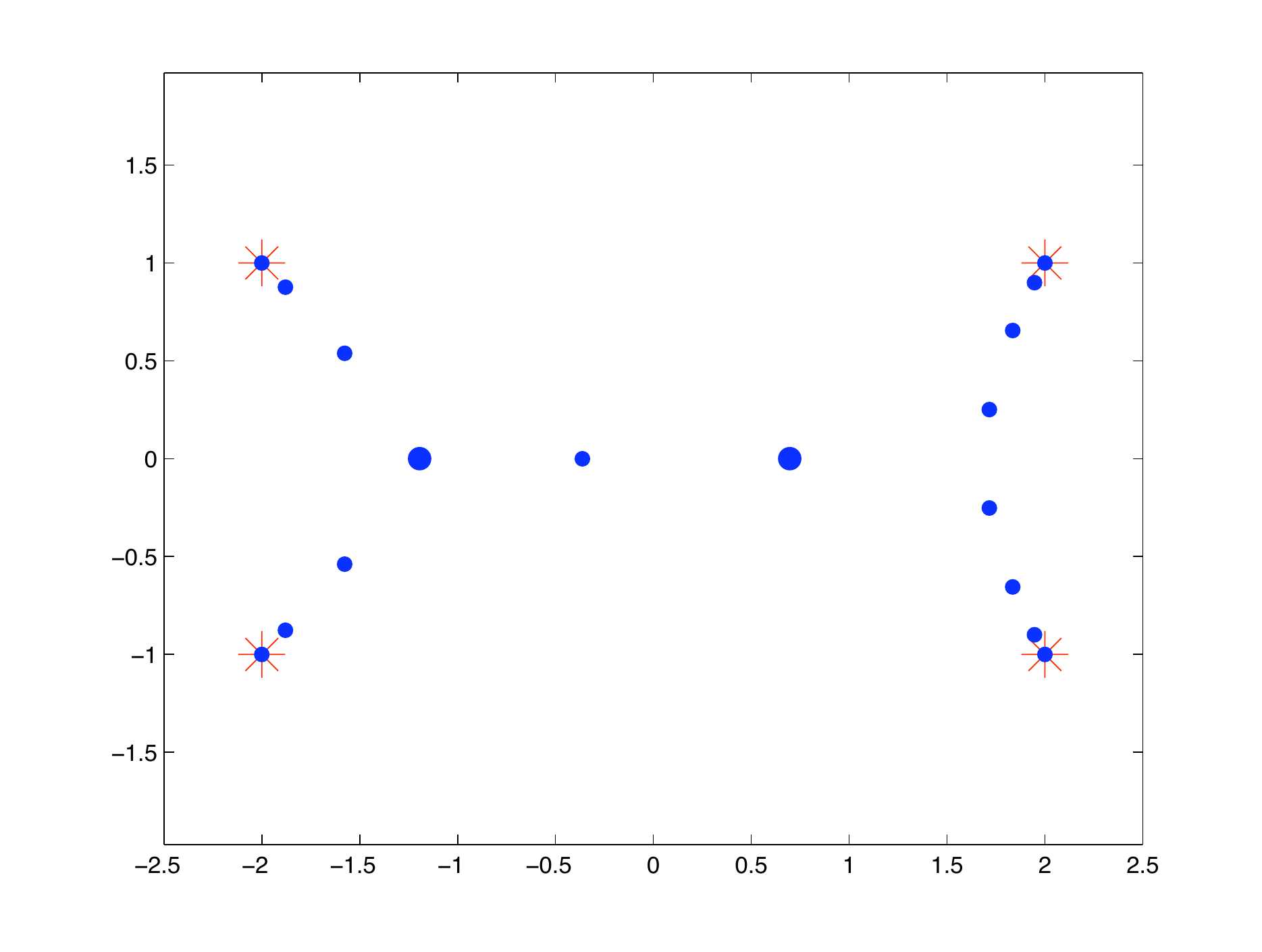}%
\end{overpic}} 
\end{tabular}
\caption{Zeros of Heine-Stieltjes (small dots) and of the corresponding Van Vleck polynomials (fat dots) for  $a_{j}$'s at vertices of a rectangle ($p=3$).}
\label{fig:experiments}
\end{figure}

\begin{figure}[htb]
\centering \mbox{\begin{overpic}[scale=0.65]%
{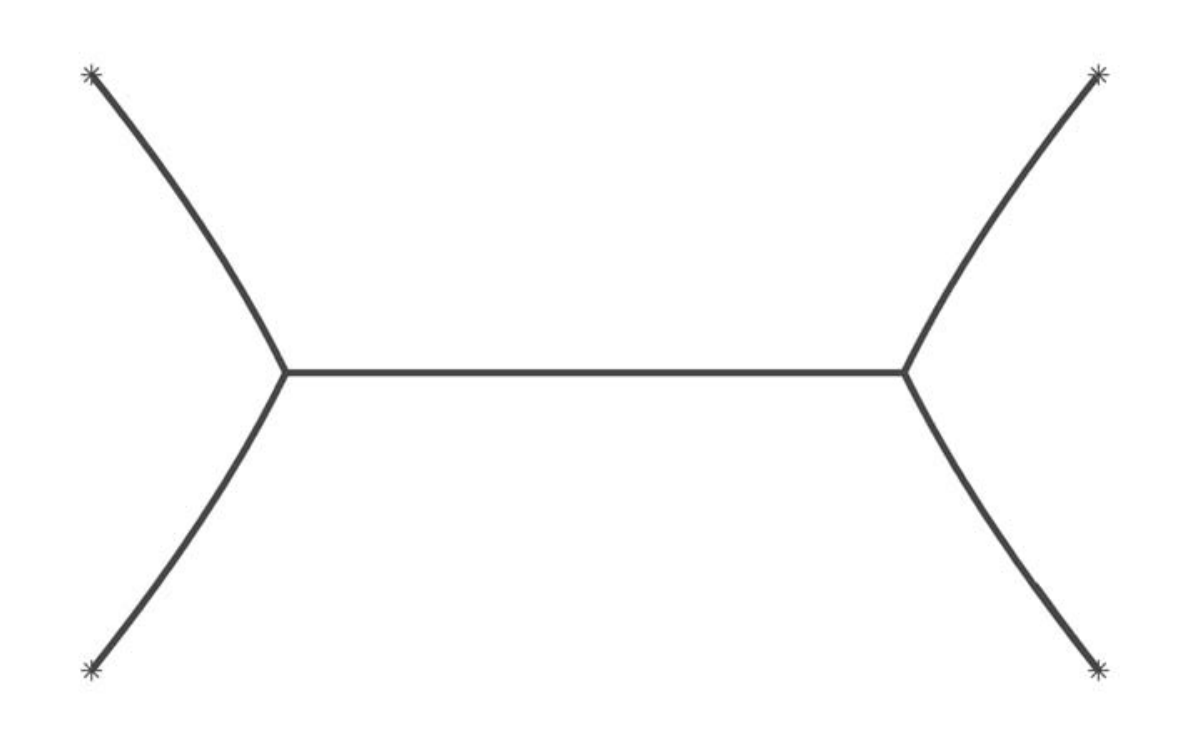}%
 \put(24,34){$v^*_1 $}
   \put(71,34){$v^*_2 $}
\put(2,55){$a_1 $}
\put(2,4){$a_2 $}
\put(94,55){$a_0 $}
\put(94,4){$a_3 $}
\end{overpic}}
 \caption{Chebotarev's compact corresponding to 4 points forming a rectangle. } \label{fig:Cheb4pts5arcsbis}
\end{figure}

Here we illustrate the situation considering the points $a_{j}$ from  $\AA$ at vertices of a rectangle whose height (vertical size) is smaller then length (horizontal size), see Fig.~\ref{fig:experiments}, with the results of some numerical experiments, as well as Fig.~\ref{fig:Cheb4pts5arcsbis}, where the corresponding Chebotarev's continuum is depicted.

In this case $V_{n}(z)=(z-v_{1,n})(z-v_{2,n})$; its zeros are determined completely by a system of two independent equations of the form \eqref{VVequationsGeneral}.

We claim that in this situation there are exactly $8+1$ homotopically different groups of systems, see Figure~\ref{fig:cases}. In eight of them we integrate along two arcs $\gamma_{k,n}$, each connecting a point from $\AA$ with a zero of $V_{n}$. In the remaining case (Figure~\ref{fig:cases}, bottom right) one of the curves will connect both zeros of $V_{n}$.

\begin{figure}[htb]
\centering \begin{tabular}{c@{\hspace{15mm}}c} 
\mbox{\begin{overpic}[scale=0.6]%
{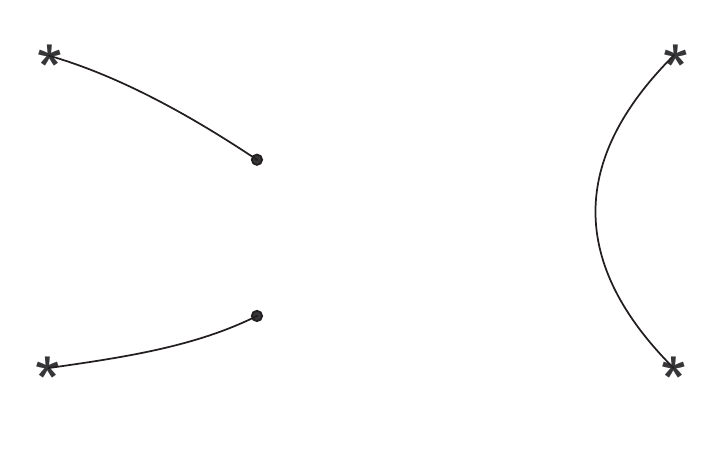}%
\put(-2.5,55){\small $a_1 $}
\put(-2.5,11){\small $a_2 $}
\put(96,11){\small $a_3 $}
\put(96,55){\small $a_0 $}
\put(38,19){\small $v_2 $}
\put(38,39){\small $v_1 $}
\end{overpic}} &
\hspace{0.5cm}
\mbox{\begin{overpic}[scale=0.6]%
{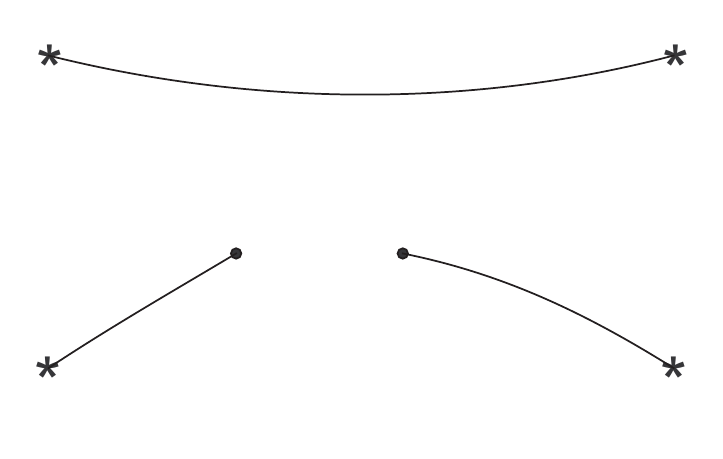}%
\put(-2.5,55){\small $a_1 $}
\put(-2.5,11){\small $a_2 $}
\put(96,11){\small $a_3 $}
\put(96,55){\small $a_0 $}
\put(47,30){\small $v_2 $}
\put(35,30){\small $v_1 $}\end{overpic}}
\\[2mm] 
\mbox{\begin{overpic}[scale=0.6]%
{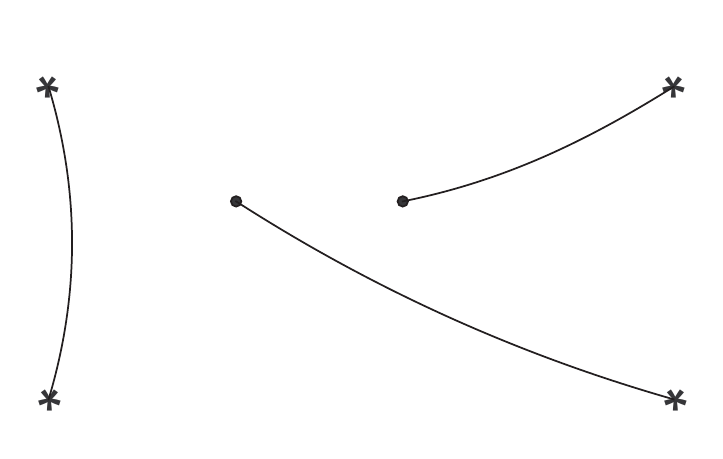}%
\put(-2.5,55){\small $a_1 $}
\put(-2.5,11){\small $a_2 $}
\put(96,11){\small $a_3 $}
\put(96,55){\small $a_0 $}
\put(48,39){\small $v_2 $}
\put(30,39){\small $v_1 $}
\end{overpic}} &
\hspace{0.5cm}
\mbox{\begin{overpic}[scale=0.6]%
{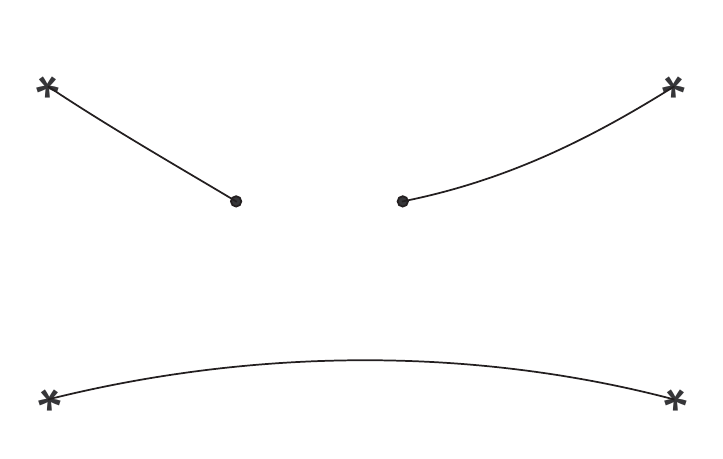}%
\put(-2.5,55){\small $a_1 $}
\put(-2.5,11){\small $a_2 $}
\put(96,11){\small $a_3 $}
\put(96,55){\small $a_0 $}
\put(48,39){\small $v_2 $}
\put(30,39){\small $v_1 $}\end{overpic}} \\[2mm] 
\mbox{\begin{overpic}[scale=0.6]%
{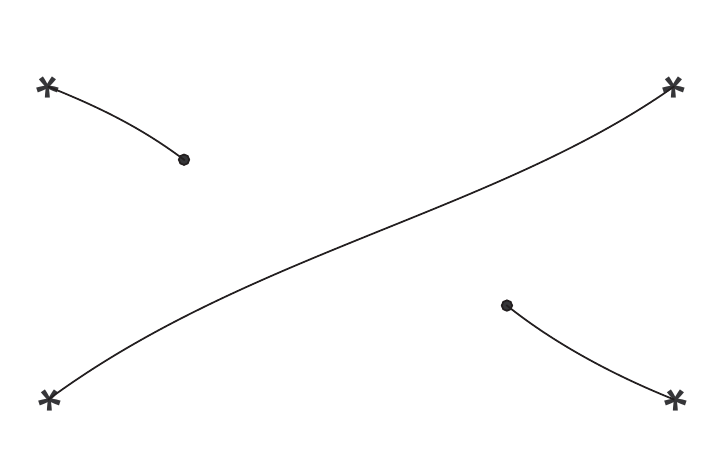}%
\put(-2.5,55){\small $a_1 $}
\put(-2.5,11){\small $a_2 $}
\put(96,11){\small $a_3 $}
\put(96,55){\small $a_0 $}
\put(62,20){\small $v_2 $}
\put(28,39){\small $v_1 $}
\end{overpic}}
 &
\hspace{0.5cm}
\mbox{\begin{overpic}[scale=0.6]%
{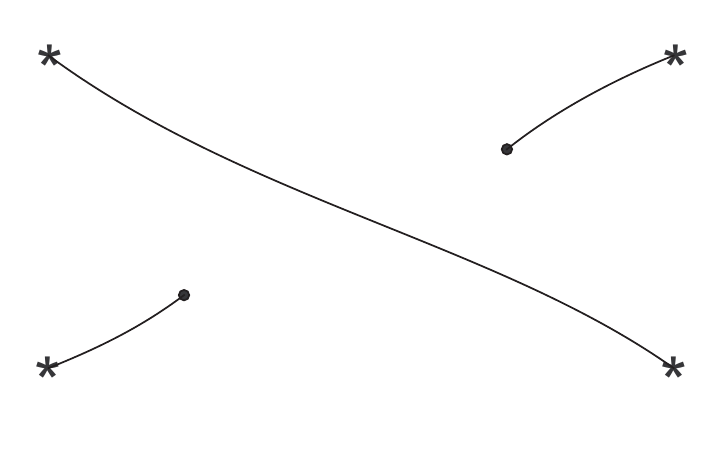}%
\put(-2.5,55){\small $a_1 $}
\put(-2.5,11){\small $a_2 $}
\put(96,11){\small $a_3 $}
\put(96,55){\small $a_0 $}
\put(62,39){\small $v_2 $}
\put(28,20){\small $v_1 $}
\end{overpic}}
 \\[2mm] 
\mbox{\begin{overpic}[scale=0.6]%
{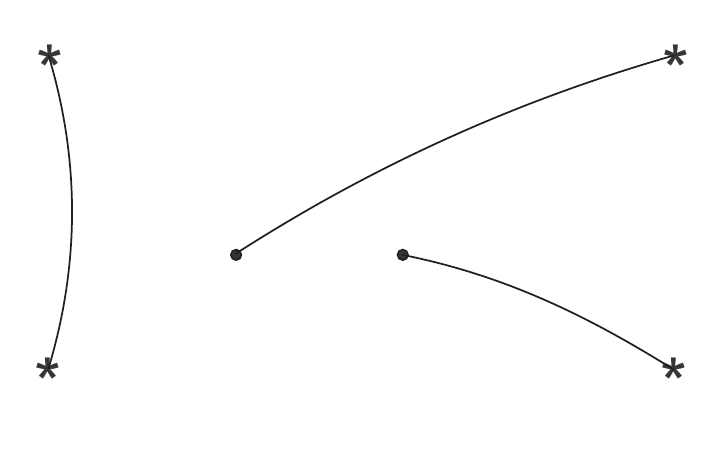}%
\put(-2.5,55){\small $a_1 $}
\put(-2.5,11){\small $a_2 $}
\put(96,11){\small $a_3 $}
\put(96,55){\small $a_0 $}
\put(49,23){\small $v_2 $}
\put(28,23){\small $v_1 $}
\end{overpic}} &
\hspace{0.5cm}
\mbox{\begin{overpic}[scale=0.6]%
{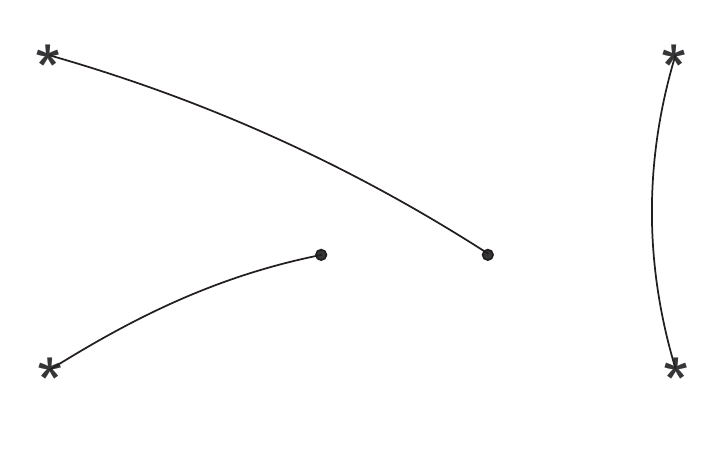}%
\put(-2.5,55){\small $a_1 $}
\put(-2.5,11){\small $a_2 $}
\put(96,11){\small $a_3 $}
\put(96,55){\small $a_0 $}
\put(63,23){\small $v_2 $}
\put(42,23){\small $v_1 $}
\end{overpic}} \\[2mm] 
\multicolumn{2}{c}{
\mbox{\begin{overpic}[scale=0.6]%
{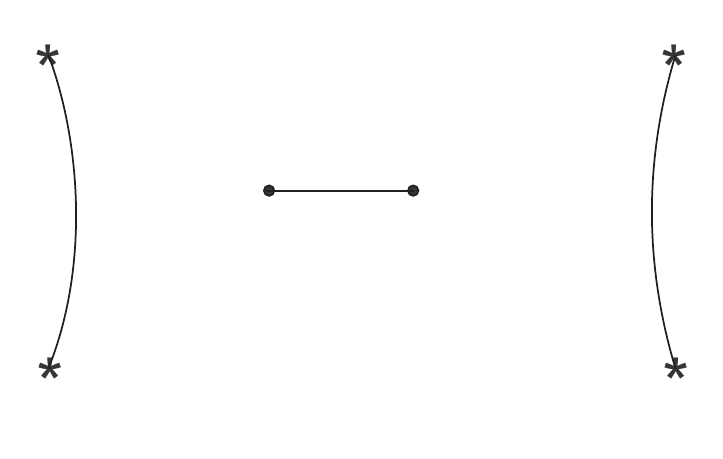}%
\put(-2.5,55){\small $a_1 $}
\put(-2.5,11){\small $a_2 $}
\put(96,11){\small $a_3 $}
\put(96,55){\small $a_0 $}
\put(55,30){\small $v_2 $}
\put(34,30){\small $v_1 $}
\end{overpic}}}
\end{tabular}
\caption{Possible homotopic classes of curves $\gamma_{n,k}$ for the case $p=3$; compare with the examples depicted in Figure~\ref{fig:experiments}.}
\label{fig:cases}
\end{figure}

For the case depicted in Figure~\ref{fig:cases}, upper left, system \eqref{VVequationsGeneral} is written as
$$
\begin{cases}\displaystyle
w_{1}(c_{1}, c_{2})=\frac{1}{\pi i}\, \int_{a_{1}}^{v_{1}} \sqrt{\frac{V_{n}(t)}{A(t)}}\, dt =
\frac{  m_{k,1} }{\lambda_{n}}+\frac{1}{\lambda_{n}}\,  \frac{B(a_{1})}{2 A'(a_{1})}+  \frac{     \delta_{k,n,1}}{\lambda_{n}},& \\[3mm] \displaystyle
w_{2}(c_{1}, c_{2})=\frac{1}{\pi i}\, \int_{a_{2}}^{v_{2}} \sqrt{\frac{V_{n}(t)}{A(t)}}\, dt =\frac{  m_{j,2} }{\lambda_{n}}+\frac{1}{\lambda_{n}}\,  \frac{B(a_{2})}{2 A'(a_{2})}+  \frac{     \delta_{j,n,2}}{\lambda_{n}}, &
\end{cases}
$$
where $(m_{k,1}, m_{j,2})$ is a pair of integers. We postpone the discussion about the possible range of variation of these constants   until next section. In order to claim that
$$
|\delta_{k,n,j}|\leq \frac{C}{n}
$$
we need to leave aside again a fraction of zeros of Van Vleck  polynomials (for which conditions of Theorem \ref{mainThm} are not satisfied), which is typically $o(n^{2})$ (recall that the cardinality of the set of all Van Vleck  polynomials corresponding to a Heine-Stieltjes polynomial of degree $n$ is $\sigma(n)=\mathcal O( n^{2})$).

\section{The Van Vleck set (outline of the proof)}

This section is essentially based on the results of our recent paper  \cite{Martinez-Finkelshtein:2009ad}, and we will be regularly referring the reader to some parts of this work for further details.

We have seen that even in the simplest cases the geometry behind the set of equations \eqref{VVequationsGeneral} is quite involved; furthermore, we have not clarified yet the selection of the index set $\mathbb M_n$. However, the situation becomes much more clear when we take limit as $n\to \infty$.

We start by considering an arbitrary $p$ for the price of omitting some non-essential details; a more thorough analysis will be carried out at the end for the particular (but non-trivial) case $p=2$.

\subsection{Zero distribution of Heine-Stieltjes polynomials and critical measures}

It is known (see e.g. \cite{Shapiro2008a}) that  the zeros of the Heine-Stieltjes polynomials accumulate on the convex hull of $\AA$. Hence, without loss of generality we may assume that there is
a subsequence $\Lambda \subset \N$ such that the Van Vleck
polynomials have a limit,
\begin{equation}\label{limitV}
\lim_{n\in \Lambda } V_n (z)=V(z)=\prod_{j=1}^{p-1} (z-v_j)\,.
\end{equation}
This limit induces the following rational quadratic differential on the Riemann sphere,
\begin{equation}\label{defQDLimit}
\varpi = -\frac{V(z)}{A(z)}\, (dz)^2.
\end{equation}

Let us consider the corresponding sequence of Heine-Stieltjes polynomials $Q_n$, $n\in \Lambda$; from Proposition \ref{prop:neclaceHS} it follows that the zeros of $Q_n$'s lie asymptotically on critical trajectories of $\varpi$, and that these trajectories have the total $\varpi$-length 1. Further results can be obtained using the electrostatic interpretation of these zeros; it has been proved in \cite{Martinez-Finkelshtein:2009ad} that any weak-* limit of the normalized zero-counting measure for $Q_n$'s is a continuous critical measure with respect to the set $\AA$ of fixed points on the plane. Let us recall the definition in its simplest form, sufficient for what follows.

With every (real-valued) Borel measure $\mu$ on $\C$  we associate its (continuous) logarithmic energy
\begin{equation}\label{defEnergyContinuous}
    E(\mu) \isdef  \iint \log \frac{1}{|x-y|}\, d\mu(x) d\mu(y)\,.
\end{equation}
Any smooth complex-valued function $h$ in the closure  $\overline \Omega $ of a domain $\Omega$ containing $\AA$ generates a local variation of $\Omega$ by $ z
\mapsto z^t=z+ t \, h(z)$, $t\in \C$, and consequently, a variation of sets $ e \mapsto
e^t \isdef \{z^t:\, z \in e\}$, and (signed) measures: $ \mu \mapsto
\mu^t$, defined by $\mu^t(e^t)=\mu(e) $.

We say that a signed measure $\mu$ supported in $ \Omega $ is a continuous \emph{$ \mathcal A $-critical} if for any $h$ smooth in $\Omega\setminus \AA$ such that $h\big|_\AA \equiv 0$,
\begin{equation}\label{derivativeEnergy}
    \frac{d}{dt}\, E (\mu^t)\big|_{t=0} = \lim_{t\to 0} \frac{E (\mu^t)- E (\mu)}{t}=0.
\end{equation}
$\AA$-critical measures appear quite frequently in many problems of approximation theory; any such a measure  $\mu$ can be characterized in terms of its logarithmic potential
$$
U^\mu(z)  \isdef  \int \log\frac{1}{|z-t|}\, d\mu(t)
$$
as follows:
\begin{lemma}[\cite{Martinez-Finkelshtein:2009ad}, Section 5.3]
\label{lemma:Sconditions}
The logarithmic potential of an $\AA$-critical measure $\mu$ satisfies the following properties:
\begin{enumerate}
\item[(i)] if $\supp(\mu )=\gamma_1\cup \dots\cup \gamma_s$, where $\gamma_j$ are the connected components of $\supp(\mu )$, then
    $$
    U^\mu(z)  = w_j=\const, \quad z\in \gamma_j, \quad j=1, \dots, s.
    $$
\item[(ii)] at any regular point $z\in \supp(\mu )$ (that is, such that locally at $z$, $\supp(\mu )$ is a smooth Jordan arc),
\begin{equation}\label{Sprop}
    \frac{\partial   U^\mu }{\partial n_+}\, \left(z \right) = \frac{\partial  U^\mu }{\partial n_-}\, \left(z \right),
\end{equation}
where $n_\pm$ are the normal vectors to $\supp(\mu )$ at $z$ pointing in the opposite directions.

Additionally, if $z\in \supp(\mu )\setminus \AA$ is not regular, then
\begin{equation}\label{zeroatnotreg}
    \grad   U^\mu(z)  =0.
\end{equation}
\end{enumerate}
Reciprocally, assume that a finite real measure $\mu$, whose support $\supp(\mu)$ consists of a union of a finite set of analytic arcs, $\supp(\mu)=\gamma_1\cup \dots\cup \gamma_s$, satisfies conditions (i) and (ii) above. Then $\mu$ is $ \AA $-critical.
\end{lemma}
In other words, (i) says that $\AA$-critical measures are in fact equilibrium measures in a piece-wise constant external field, exhibiting additionally (see (ii)) the so-called $S$-property, introduced first by Stahl \cite{MR88d:30004a} and in a more general context, by Gonchar and Rakhmanov \cite{Gonchar:84} (where it was used to establish the well-known ``$1/9$''-conjecture in approximation theory). A first rigorous proof of the connection of a critical measure with the $S$-property appeared in \cite{Rakhmanov94}.

With a polynomial $Q(z)=\prod_{k=1}^n (z-\zeta_{k})$ of degree $n$ we associate its (normalized) zero counting measure
$$
\nu(Q)=\frac{1}{n}\, \sum_{k=1}^{n} \delta_{\zeta_{k}}.
$$
The differential equation \eqref{DifEq} is an expression of the fact that the zeros of $Q_n$ sit in equilibrium (zeros of the gradient of the discrete total energy) in presence of an external field depending from the residues $\rho _k$ in \eqref{BoverA}. The intensity of this field decays with $n$ proportionally to $1/n$, so that in the limit we get
\begin{proposition}[\cite{Martinez-Finkelshtein:2009ad}, Section 7]
\label{the:weakAsymptoticsWithField}
Let $\nu_n=\nu(Q_{n})$ be a zero-counting measure corresponding to a sequence of Heine-Stieltjes polynomials $Q_n$. Then any weak-* limit point $\mu$ of $\nu_n$ is a unit continuous $\AA$-critical measure.
\end{proposition}
In other words, weak-* limits of the normalized zero counting measures of $Q_n$'s are unit positive $\AA$-critical measures.
The inverse inclusion (that any unit positive $\AA$-critical measure is a weak-* limit of the normalized zero counting measures of Heine-Stieltjes polynomials) is also valid, but it cannot be established within the framework of this paper. 

Continuous critical measures can be also related with rational quadratic differentials on $\overline \C$ by means of a variational argument. Our WKB analysis allows to establish this link, at least partially, more directly:
\begin{theorem}
\label{the:weakAsymptotics}
Let $y=Q_{n}$ be a Heine-Stieltjes polynomial (solution of \eqref{DifEq}) corresponding to $V_{n}$, and assume that \eqref{limitV} holds. Then the quadratic differential
$$
\varpi = -\frac{V}{A}(z) \, dz^2
$$
is closed (all its trajectories are either critical or closed), and there is a probability $\AA$-critical measure  $\mu$ on $\C$ such that the normalized zero counting measures $\nu(Q_n)$ converge (in a weak-* sense and along the subset $\Lambda$) to  $\mu$. The support $\Gamma =\supp(\mu )$ consists of critical trajectories of $\varpi$, $\C\setminus \Gamma$ is connected, and we can fix the single valued branch of $\sqrt{V/A}$ there by $\lim_{z\to\infty} z \sqrt{V(z)/  A(z) }=1$. With this convention,
\begin{equation}\label{asymptNthRoot}
    \lim_n \left| Q_n(z)\right|^{1/n} =\exp\left( \Re \int^z \sqrt{\frac{V}{A}}(t) \, dt \right)
\end{equation}
locally uniformly in $\C\setminus \Gamma$, where a proper normalization of the integral in the right hand side is chosen, so that
$$
\lim_{z\to\infty} \left( \Re \int^z \sqrt{\frac{V}{A}}(t) \, dt- \log |z|\right)=0.
$$
\end{theorem}

The analysis above yields the following addendum, that we state using the notation and assumptions of Theorems \ref{mainThm} and \ref{the:weakAsymptotics}:
\begin{proposition}
\label{mainThmLimitCase}
The support  $\Gamma =\supp(\mu )$ is comprised of $p$ analytic arcs $\gamma_{k}$, such that
	$$
	\lim_{n\in \Lambda} \Gamma_{n}=\Gamma.
	$$
\end{proposition}

\subsection{Structure of the family of positive critical measures}

Theorem \ref{the:weakAsymptotics} is essential for understanding the asymptotics of Heine-Stieltjes polynomials, although this result is in a certain sense implicit, since it depends on the limit $V$ of the Van Vleck polynomials $V_n$, that constitute therefore the main parameters of the problem. We must complement this description with the study of the set of all possible limits $V$.

To this end it is convenient to consider a correspondence between closed quadratic differentials and $\AA$-critical measures, independently of their origin. In general, this is not a one-to-one correspondence, since many critical measures may correspond to the same closed quadratic differential. It was pointed out in \cite{Martinez-Finkelshtein:2009ad} that the bijection between closed quadratic differentials of the form \eqref{defQDLimit} and signed $\AA$-critical measures is restored if we restrict ourselves to signed measure with a connected complement of the support. Indeed, any such a measure is supported on a finite union of analytic arcs, and these arcs are necessarily critical trajectories of a quadratic differential like in \eqref{defQDLimit}. Reciprocally, from \eqref{defQDLimit} we can construct explicitly (using the Sokhotsky-Plemelj formulas) a measure $\mu$ that satisfies \eqref{derivativeEnergy}.

Thus, any closed quadratic differential uniquely generates an $\AA$-critical measure $\mu$ with a connected complement to its support. In general, such a $\mu$ is a signed measure, and for our purposes one has to select those quadratic differentials associated with positive measures. One of the ways to solve this problem goes through describing the global structure of the trajectories of closed rational quadratic differentials with fixed denominators on the Riemann sphere, and the corresponding parameters (numerators), and later extracting from this set the subset giving rise to positive unit $\AA$-critical measures.

This way is implemented in for $p=2$ at the end of this paper. For an arbitrary $p$ such an investigation would be more difficult to carry out, in particular due to the heavy combinatorics involved. Below we present a direct construction from \cite{Martinez-Finkelshtein:2009ad}, based on the $v$-local coordinates that we introduce in the space of the $\AA$-critical measures, as well as on the topological structure of the Chebotarev's compact for $\AA$. We start by introducing the $v$-coordinates.

Let us recall the notation. We have the fixed set $\AA=\{a_0, a_1, \dots, a_p\}$ of distinct points on $\C$, $A(z)=\prod_{j=0}^p (z-a_j)$. For any signed $\AA$-critical measure $\mu$ there exists a rational quadratic differential $\varpi$ on the Riemann sphere $\overline \C$ given by
\begin{equation}\label{quadDiffBis}
\varpi  (z)=- R(z)\, (dz)^2, \qquad R(z)\isdef \frac{V(z)}{  A(z) },  \qquad V(z)\isdef \prod_{j=1}^{p-1} (z-v_j),
\end{equation}
such that $\supp (\mu)=\Gamma_{\mu}=\Gamma=\gamma _1 \cup \dots \cup \gamma _p$ is a union of trajectories of $\varpi$.

We begin by introducing the local coordinates under the assumption that measures are in general position; this
notion of genericity (as opposed to some more special or coincidental cases that are possible) means in our context that $\Gamma$ is comprised of exactly $p$ disjoint arcs. In consequence, zeros $v_j$'s of $V$ are simple, and vector $\vt{v}\isdef \{v_1, \dots, v_{p-1} \}\in \C^{p-1}$  can be used as a local coordinate. Define also
$$
\mathcal V \isdef \{ \vt{v}:\, \varpi \text{ is closed}\}.
$$

Next we use these local coordinates to introduce the period mapping for our quadratic differentials. The Carath\'{e}odory boundary of $\overline \C\setminus \Gamma $ consists of $p$ components $\widehat{\gamma}_k \isdef \gamma _k^+ \cup \gamma _k^-$, with a positive orientation with respect to $\overline \C\setminus \Gamma $. We can consider $\widehat{\gamma}_k$ as cycles in $\overline{\C}\setminus \Gamma$ enclosing the endpoints of $\gamma _k$. Part of $\mathcal R$ over $\overline{\C}\setminus \Gamma$ splits into two disjoint sheets, so we may consider $\widehat{\gamma}_k$ as cycles on $\mathcal R$.

Let us define
\begin{equation}\label{defW_k9}
    w_k(\vt v)=w_k (\vt v, \Gamma ) \isdef \frac{1}{2\pi i}\, \oint_{\widehat{\gamma}_k} \sqrt{R(z)}dz, \quad k=1, \dots, p,
\end{equation}
where $\sqrt{R}\big|_{\widehat{\gamma}_k}$ are the boundary values of the branch of $\sqrt{R}$ in $\overline{\C}\setminus \Gamma $ defined by $\lim_{z\to \infty} z \sqrt{R(z)} =1$. Clearly, the boundary values $(\sqrt{R})_\pm$ on $\gamma_k^\pm$ are opposite in sign. Therefore, with any choice of orientation of $\gamma _k$ and a proper choice of $\sqrt{R}=(\sqrt{R})_+$ on $\gamma_k$, we will have
\begin{equation}\label{periods9}
    w_k(\vt v)=\frac{1}{\pi i}\, \int_{\gamma _k} \sqrt{R(z)} dz, \quad k=1, \dots, p.
\end{equation}
By the Cauchy residue theorem we have that $w_1+\dots+w_p=1$ for any $\vt v\in \C^{p-1}$. Thus, we can restrict the mapping $\vt v\mapsto \vt w$ to $p-1$ components of $\vt w \isdef (w_1, \dots, w_{p-1})\in \C^{p-1}$. In this way, we have defined the mapping
\begin{equation}\label{mapping9}
  \mathcal P(\cdot, \Gamma ):\, \C^{p-1} \to  \C^{p-1} \quad \text{such that} \quad  \mathcal P(\vt v, \Gamma ) = \vt w(\vt v, \Gamma ) .
\end{equation}
Each component function $w_j(v_1, \dots, v_{p-1})$ is analytic in each coordinate $v_k$ (even if $v_k$ is at one of the endpoints of $\gamma _k$). Once defined by the integral in \eqref{periods9}, this analytic germ allows an analytic continuation along any curve in $\C\setminus \AA$. Arcs $\gamma _k$ are not an obstacle for the continuation since the integral in \eqref{periods9} depends only on the homotopic class of $\Gamma$  in $\C\setminus (\AA\cup \vt v)$. The homotopy of $\Gamma $ is a continuous modification of all components simultaneously in such a way that they remain disjoint in all intermediate positions. Under this assumption we can continuously modify the selected branch of $\sqrt{R}$ in $\overline{\C}\setminus \Gamma $ along with the motion of $\Gamma $.

\begin{proposition}
\label{prop:9_1}
Mapping $w=\mathcal P(\vt v)$ is locally invertible at any $\vt v(v_1, \dots, v_{p-1})\in (\C\setminus \AA)^{p-1}$ with $v_i\neq v_j$ for $i\neq j$.
\end{proposition}
Furthermore, we have
\begin{proposition}
\label{prop:9_2}
Let $\mu^0$ be an $\AA$-critical measure such that $\Gamma=\supp(\mu ^0)$ has connected components $\gamma _1^0,  \dots , \gamma _{p}^0$, and $\overline{\C}\setminus \Gamma$ is connected. Let $\varpi_0 =R_0(z)(dz)^2$ be the quadratic differential associated with $\mu^0$, where $R_0=V_0/A$, and $\vt v^0=(v_1^0, \dots, v_{p-1}^0)$ is the vector of zeros of $V_0$. Assume that $\varpi_0$ is in general position (that is, all $v_k$'s are pairwise distinct and disjoint with $\AA$). Then for an $\varepsilon >0$ and any $m_j \in\R$, $j\in \{1, \dots, p-1\}$, satisfying
$$
|m_j - \mu ^0(\gamma _j^0)|<\varepsilon , \quad j=1, \dots, p-1,
$$
there exists a unique solution $\vt v \in \mathcal V$ of the system
\begin{equation}\label{systemforMu9}
    w_j(\vt v,\Gamma )=m_j, \quad j=1, \dots, p-1.
\end{equation}
The quadratic differential $\varpi = -R(z)(dz)^2$, $R(z)=\prod_{k=1}^{p-1} (z-v_k)/A(z)$, is closed, and the associated $\AA$-critical measure $\mu $ satisfies $\mu (\gamma _j)=m_j$, $j=1, \dots, p-1$, where $\supp(\mu )=\gamma _1\cup \dots \cup \gamma _{p }$ and $\supp(\mu)$ is homotopic to $\supp(\mu^0 )$.
\end{proposition}

This Proposition introduces a topology in the set of $\AA$-critical measures in general position. We call a \emph{cell} any connected component of this topological space. A measure in general position preserves sign along any connected component of its support. This sign is subsequently preserved when we homotopically modify the measure within its cell. In particular, if $\mu$ is a positive $\AA$-critical measure, then all measures in the same cell with $\mu$ are positive.

In $v$-coordinates a cell $G=G(\Gamma)$ in $\mathcal V$ is a subspace of $\C^{p-1} \simeq \R^{2p-2}$, which is a manifold of the real dimension $p-1$, defined by $p-1$ real equations of the form
\begin{equation}\label{Takemura9}
    \Im w_j(v)=\Im \left(\frac{1}{2\pi i}\, \int_{\widehat \gamma _j} \sqrt{R(t)}\, dt\right) =0, \quad j=1, \dots, p-1,
\end{equation}
where $\widehat{\Gamma }=\widehat{\gamma }_1\cup \dots \cup \widehat{\gamma }_{p-1} \cup \widehat{\gamma }_p$ is a union of Jordan contours $\widehat{\gamma }_k$ on $\mathcal R$ (double arcs) depending on $v$, but mutually homotopically equivalent for values of $v$ from the same cell. Practically, any $v_0\in G$ has a neighborhood of $v$ satisfying \eqref{Takemura9} with constant $\widehat{\gamma }_k$'s.

Next, let $\mathcal V_+$ be the set of positive $\AA$-critical measures; it will be comprised of $3^{p-1}$ cells. Below we introduce these cells by means of the Chebotarev's continuum  $\Gamma ^*=\Gamma ^*(\AA)$, associated with $\AA=\{ a_0, \dots, a_p\}$ (that is, the continuum solving the extremal problem \eqref{defChebotarev}). It constitutes the ``universal boundary'', i.e., the intersection of all boundary points of all cells. In general, isolating a cell from its boundary point is a formidable task. We take advantage here that we are interested only in those that yield positive $\AA$-critical measures, that are characterized by the trivial homotopy.

The Chebotarev's continuum consists of critical trajectories of $\varpi= -R^* dz^2$, $R^*=V^*/A$. Assume that the set $\AA$ satisfies the condition that $A$ and $V^*$ do not have common zeros and $V^*(z)=\prod_{k=1}^{p-1} (z-v_k^*)$ does not have multiple zeros. Thus, the critical set $\AA^*= \AA \cup \{v^*_1, \dots, v_{p-1}^* \}$ consists of $2p$ different points, and $\Gamma ^*$ is comprised of $2p-1$ arcs, that are critical trajectories of $\varpi$. Each trajectory joins two different points from $\AA^*$.

Now we come to the procedure of selection of combinatorial (rather than homotopic) types of cells; once the combinatorial type is fixed, the homotopic one will be determined from the Chebotarev's continuum, as described next. Each zero $v^*=v^*_k$, $k=1, \dots, p-1$, is connected by component arcs of $\Gamma ^*$ with three other points, say $a_1^*, a_2^*, a_3^*\in \AA^*$. We select one of these three arcs (for definiteness, $[v^*,a_1^*]$) and join two other arcs to make a single arc $[a_2^*, a_3^*]$, bypassing $v^*$ (we think that the arc $[a_2^*,a_3^*]$ still follows the two arcs from $\Gamma ^*$, but without touching $v^*$, instead passing infinitely close to it). This procedure, carried out at each zeros $v_k^*$ of $V^*$, creates a compact set $\Gamma $, and consequently, a cell $G(\Gamma )$ of corresponding measures $\mu\in \widehat{\mathcal V}_+$.

The selection of $\Gamma $, and hence, of the cell $G(\Gamma )$, is made by choosing one of the three connections for each $v^*_k$; there are $3^{p-1}$ ways to make the choice. Any choice splits $\Gamma ^*$ into $p$ ``disjoint'' arcs $\Gamma ^*=\gamma _1\cup \dots\cup \gamma _p$; out of them we select $p-1$ arcs (to make an homology basis for $\C\setminus \Gamma ^*$) and then consider the corresponding cycles $\widehat{\gamma }_k$, as described above.

Finally, we describe the cell $G(\Gamma )$ in terms of the mapping $\mathcal P$. Let $w=w(\vt v)=\mathcal P(v, \widehat{\Gamma })$;  the cell $G(\Gamma )$ is completely defined by the system
\begin{equation}\label{6of9}
    w_j(\vt v)=\mu _j \in \R_+, \quad j=1, \dots, p-1;
\end{equation}
more precisely, there exists a domain $M(\Gamma )=\{ (\mu _1, \dots, \mu _{p-1})\in \R_+^{p-1} \}$ such that for any  point $(\mu _1, \dots, \mu _{p-1})\in M(\Gamma )$ system \eqref{6of9} has a unique solution $\vt v\in \C^{p-1}$. Moreover, the corresponding measure $\mu=\mu _v$ satisfies $\mu (\gamma _j)=\mu _j$, and $\supp(\mu )=\gamma _1\cup \dots \cup \gamma _{p}=\Gamma _v$ is homotopic to $\Gamma $.

Summarizing, a rough description of the set $\widehat{\mathcal V}_+$ of unit positive $\AA$-critical measures may be made as follows. The set $\widehat{\mathcal V}_+$ is a union of $3^{p-1}$ of closed bounded cells $\overline{G}(\Gamma )$ ($\Gamma =\gamma_1\cup \dots \cup \gamma _p$ may be selected in $3^{p-1}$ ways). The interior $G(\Gamma )$ of each cell consists of measures $\mu$ in general position with $\supp(\mu )$ homotopic to $\Gamma $. Interiors of different cells are disjoint. Chebotarev's measure $\mu ^*$ (Robin measure of $\Gamma ^*$) is the only common point of all boundaries: $\mu ^*=\bigcap_\Gamma  \partial G(\Gamma )$. A graphical description of all these cells for $p=3$, obtained from the construction just described, is contained in Figure \ref{fig:cases}.

\subsection{Case $p=2$} \label{subsec:VVforP=2}

In order to clarify  the construction above let us discuss in more detail the simplest (but far from trivial) case of $p=2$.

Let us introduce the following set of the plane. For the quadratic differential
\begin{equation}\label{qdp2}
    \varpi_v= \frac{v-z}{A(z)}\, dz^2
\end{equation}
define
$$
\mathcal V \isdef \left\{ v\in \C:\, \varpi_v \text{ is closed} \right\},
$$
as well as the \emph{Van Vleck set}
$$
\mathcal V_+ \isdef \left\{ v\in \C:\, v \text{ is an accumulation point of the zeros of Van Vleck polynomials} \right\}.
$$
A direct consequence of Theorem \ref{the:weakAsymptotics} is that
$$
\mathcal V_+ \subset \mathcal V.
$$
As it was mentioned for the general case, this inclusion is proper.

As in Subsection \ref{subs:p=2Cheb}, we denote by $v^*$ the Chebotarev's center of $\AA$ (the value of $v$ in \eqref{qdp2} such that $\varpi_{v^*}$ corresponds to the Chebotarev's continuum for $\AA$), and let $\AA^*=\AA \cup \{ v^{*}\}$.

\begin{theorem}[\cite{Martinez-Finkelshtein:2009ad}, Section 8.3] \label{cor:structureV}
The set $\mathcal V$ is a union of a countable number of analytic arcs $\ell_k$, $k\in \Z$, each connecting $v^*$ and $\infty$.

Two arcs from $\mathcal V$ are either identical or have $v^*$ as the only finite common point. The homotopic type of the critical trajectories of $\varpi_v$ in $\C\setminus  \AA $ remains invariant on each arc $\ell_k\setminus \AA^*$.

There are three distinguished arcs $\ell_k$, $k\in \{0, 1, 2\}$, such that
\begin{enumerate}
\item[(i)] $\ell_k$ connects $v^*$ with infinity and passes through $a_k$;
\item[(ii)] for every $v\in \ell_k$ the homotopic class of trajectories of the closed quadratic differential $\varpi_v$ is trivial.
\end{enumerate}
\end{theorem}
In the terminology introduced in the previous subsection, this theorem says that each cell (connected component) in the topological space of $\AA$-critical measures is homeomorphic to an analytic arc with endpoints either at $v^{*}$, $\AA$ or infinity. This settles the problem of identification of all closed quadratic differentials $\varphi_{v}$. Using the arguments described above we single out the differentials corresponding to \emph{positive} measures:
\begin{theorem} \label{thm:structureVPlus2}
Let $\ell_k$, $k\in \{0, 1, 2\}$, be the distinguished arcs in $\mathcal V$ described in Theorem  \ref{cor:structureV}. The set  $\widehat{\mathcal V}_+$ is the union of the sub-arcs $\ell_k^+$ of each $\ell_k$, $k\in \{0, 1, 2\}$, connecting $a_k$ with the Chebotarev's center $v^*$ (and lying in the convex hull of $\AA$).

Furthermore, if  $v\in \ell_k \cap \widehat{\mathcal V}_+$, $k\in \{ 0, 1, 2\}$, then there is a critical trajectory $\gamma(v)$ of $\mu_{v}$ connecting $v$ with the pole $a_k$ and such that
\begin{equation}\label{inequalityMu}
    0\leq \mu_v(\gamma(v)) \leq M_k,
\end{equation}
where $M_{k}$ is defined in \eqref{defMuK}.
In this case both critical trajectories of $\varpi_{v}$ that constitute the support of $\mu_{v}$ are homotopic to a segment. The bijection $\mu_v(\gamma(v)) \leftrightarrow v$ is a parametrization of the set $\ell_k \cap \widehat{\mathcal V}_+$ by points of the interval $[0, m _k]$.
\end{theorem}

In other words, the only three cells corresponding to positive $\AA$-critical measures are homeomorphic to three analytic arcs, joining the Chebotarev's center $v^{*}$ with each pole $a_{j}\in \AA$ (in fact, these arcs form a star homeotopic and quite close to the Chebotarev's continuum, although in general not coincident with it). If $v$ travels such an arc, the boundary of the cell is reached when either endpoint of this arc is met. If we continue further along the same arc, a new cell is entered, corresponding to  sign-changing $\AA$-critical measures (see an illustration of the correspondence between the position of $v$ on $\mathcal V$ and the trajectories of $\varpi_v$ in Figure \ref{fig:cases}).

\begin{figure}[htb]
\centering \begin{tabular}{l@{\hspace{15mm}}c} 
\mbox{\begin{overpic}[scale=0.3]%
{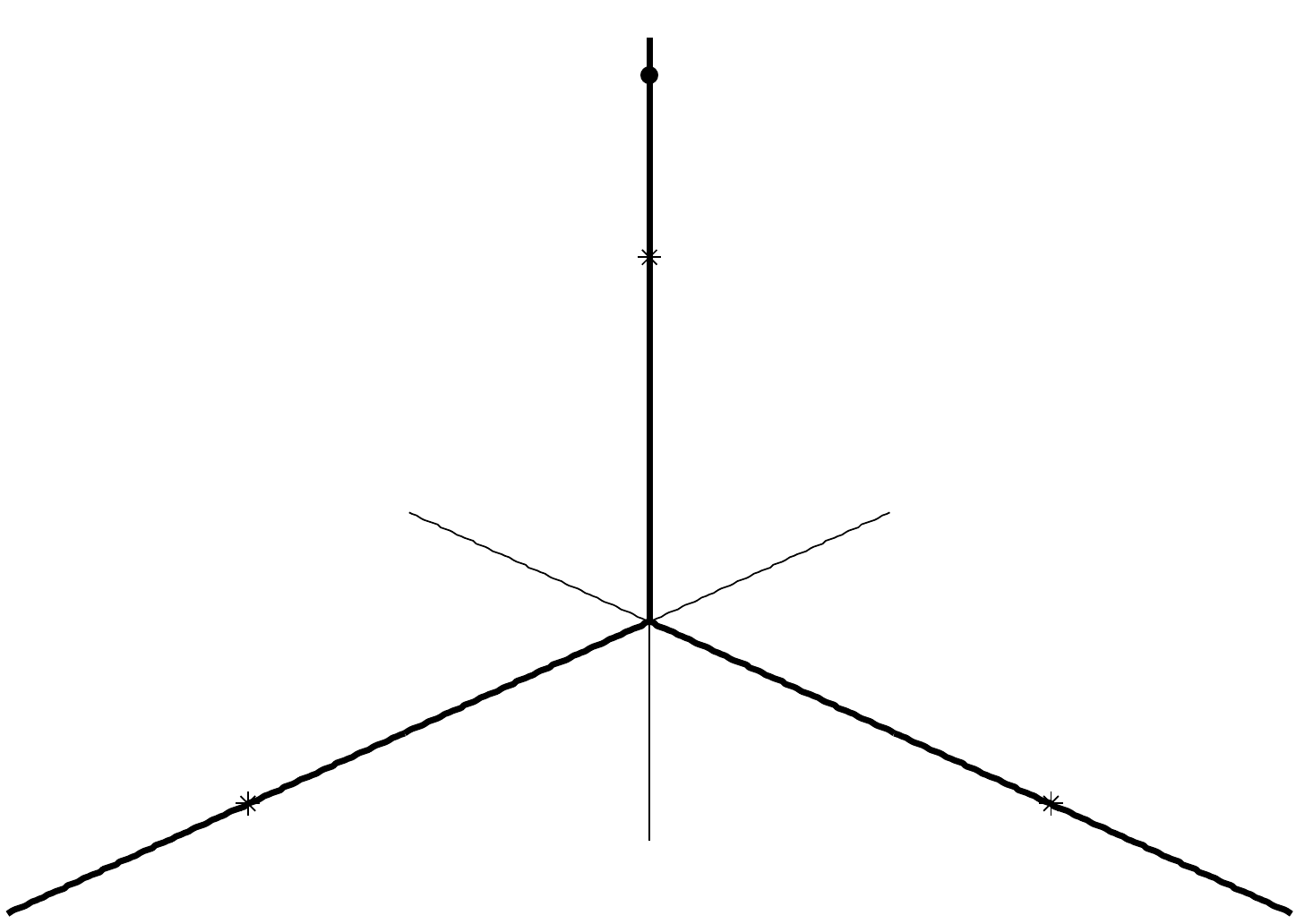}%
\put(58,45){\emph{i)} $\mu_v(\gamma(v))<0$}
      \put(40,49){\small $a_k $}
      \put(40,63){\small $v $}
\end{overpic}} &
\mbox{\begin{overpic}[scale=0.35]%
{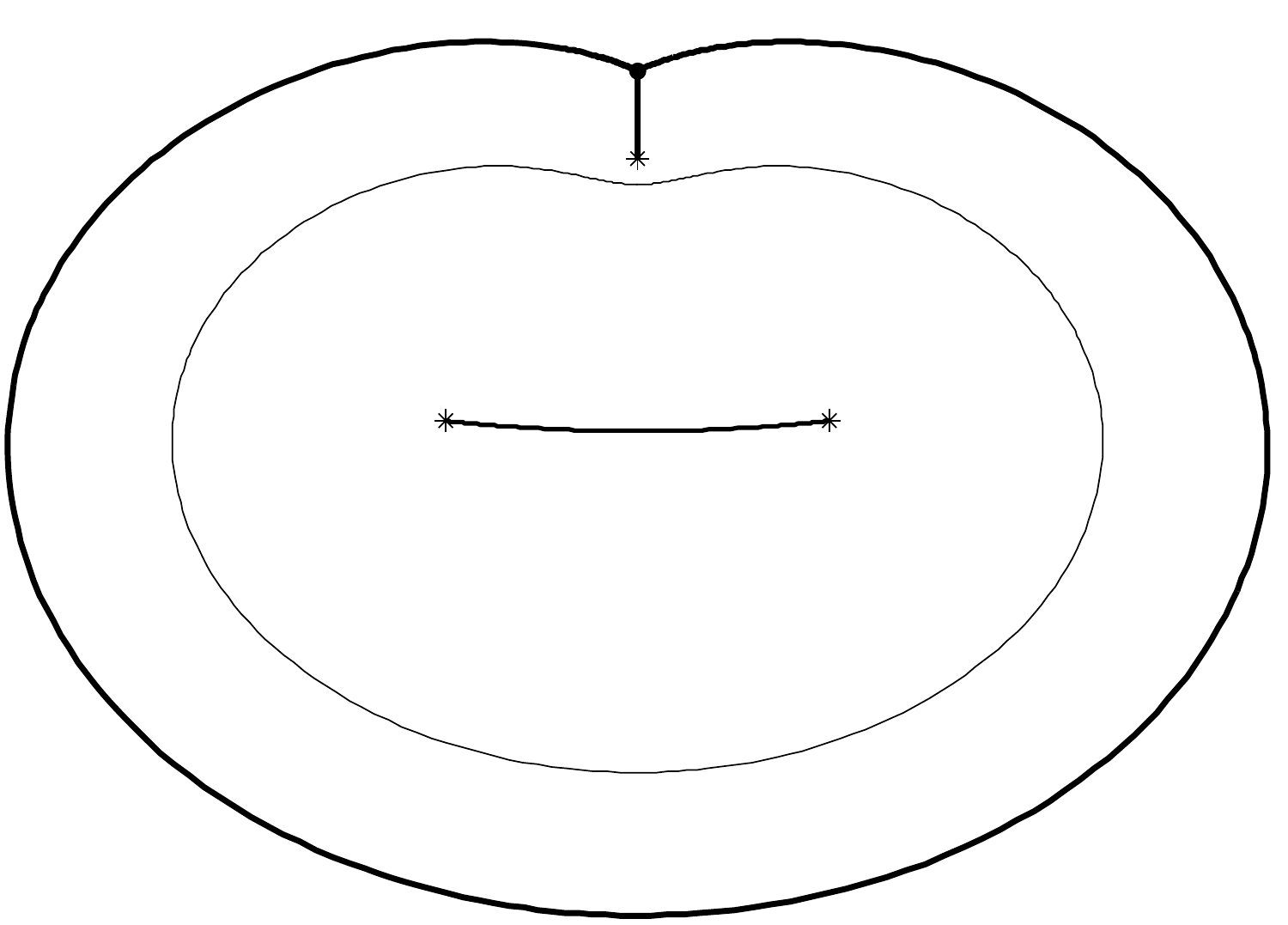}%
\end{overpic}}
\\ 
\mbox{\begin{overpic}[scale=0.3]%
{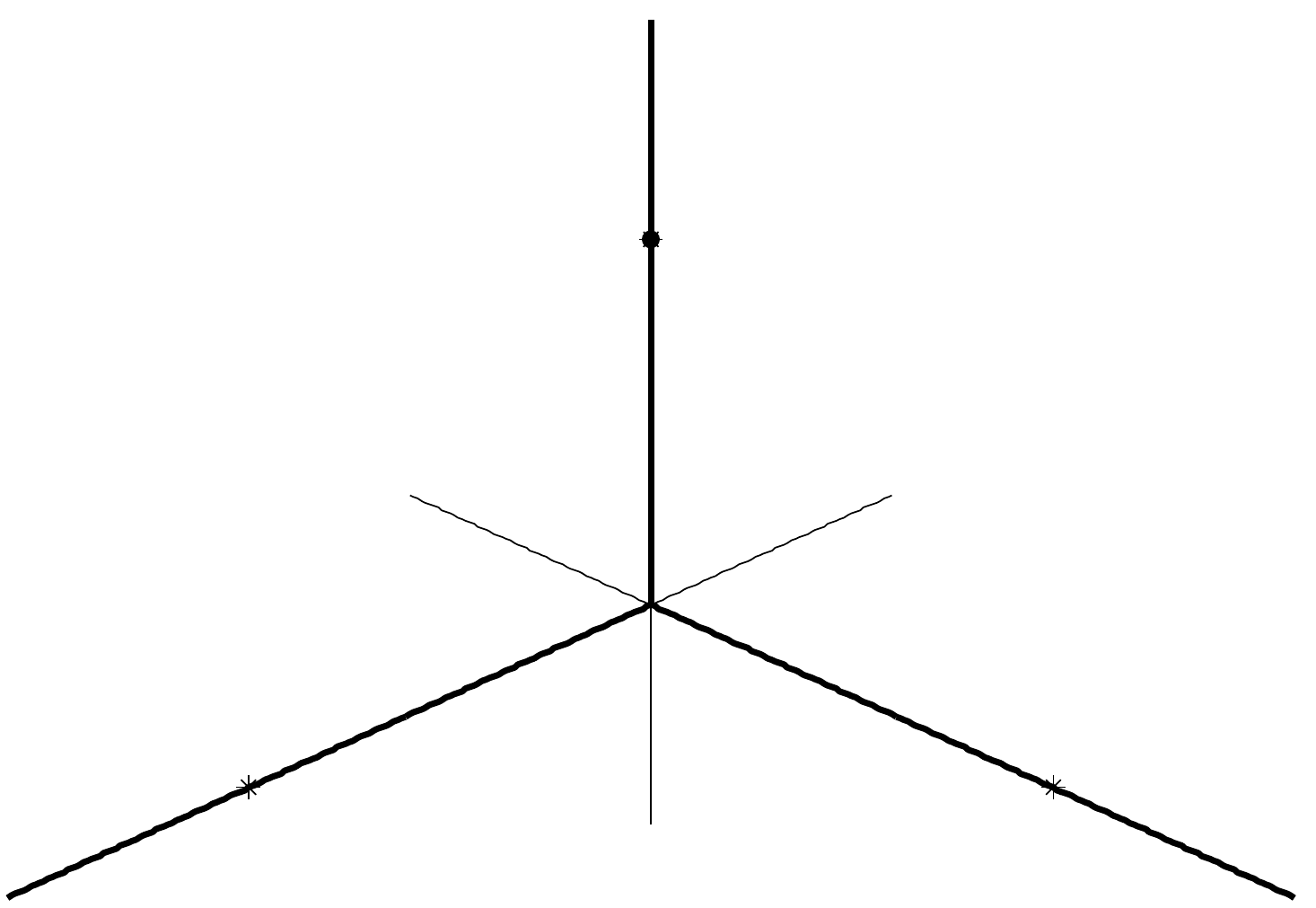}%
\put(58,45){\emph{ii)} $ \mu_v(\gamma(v))=0$}
      \put(25,49){\small $v=a_k $}
\end{overpic}} &
\mbox{\begin{overpic}[scale=0.35]%
{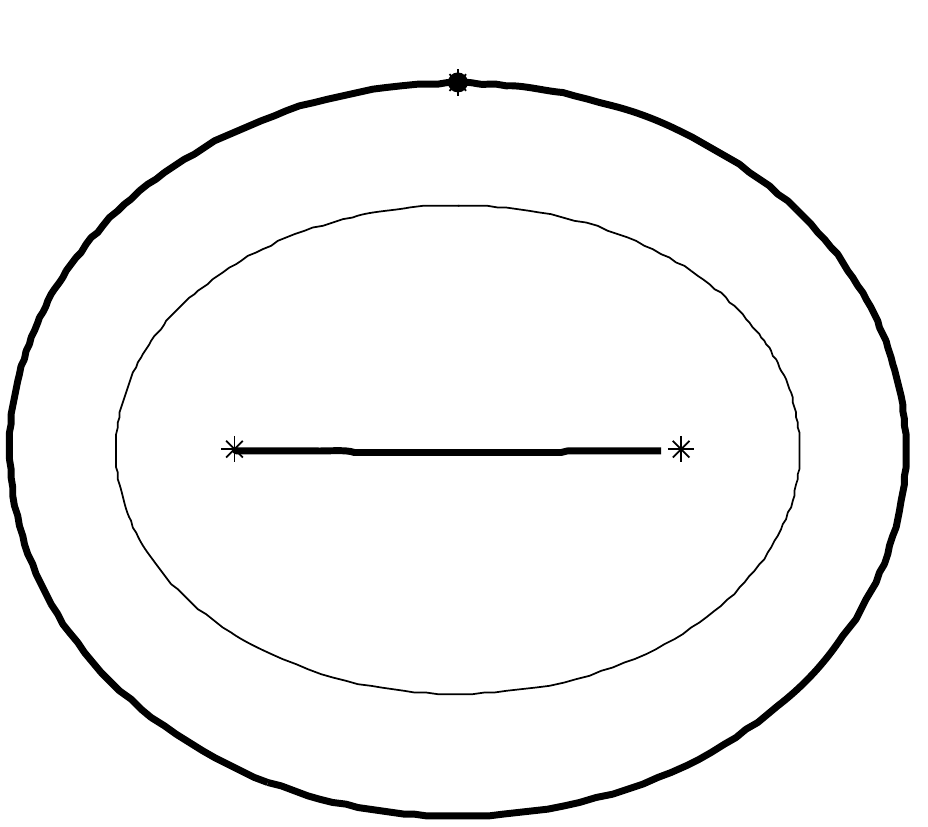}%
\end{overpic}} \\ 
\mbox{\begin{overpic}[scale=0.3]%
{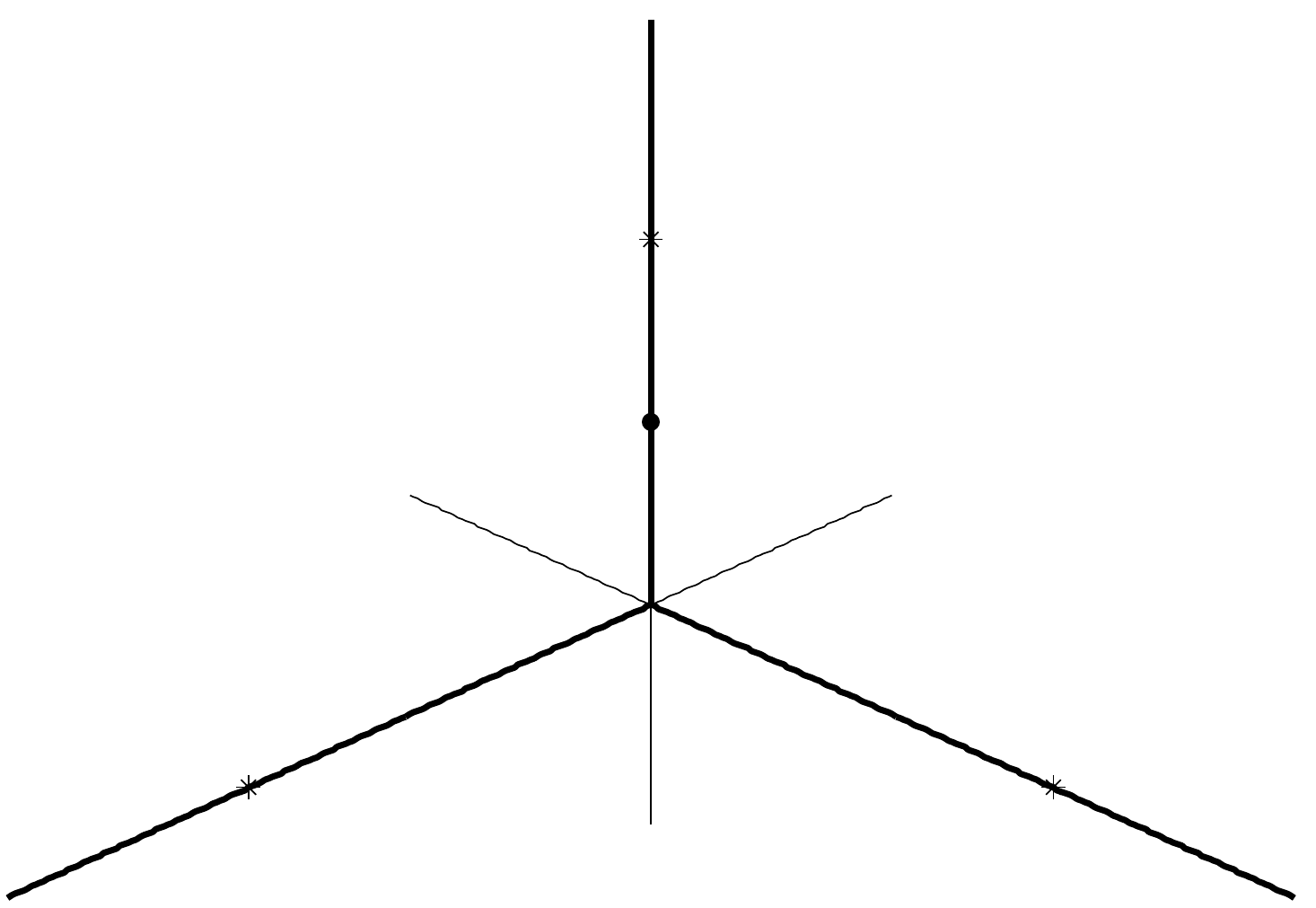}%
\put(58,45){\emph{iii)} $0<\mu_v(\gamma(v))<m_k$}
      \put(40,49){\small $a_k $}
      \put(43,35){\small $v $}
\end{overpic}}
 &
\mbox{\begin{overpic}[scale=0.35]%
{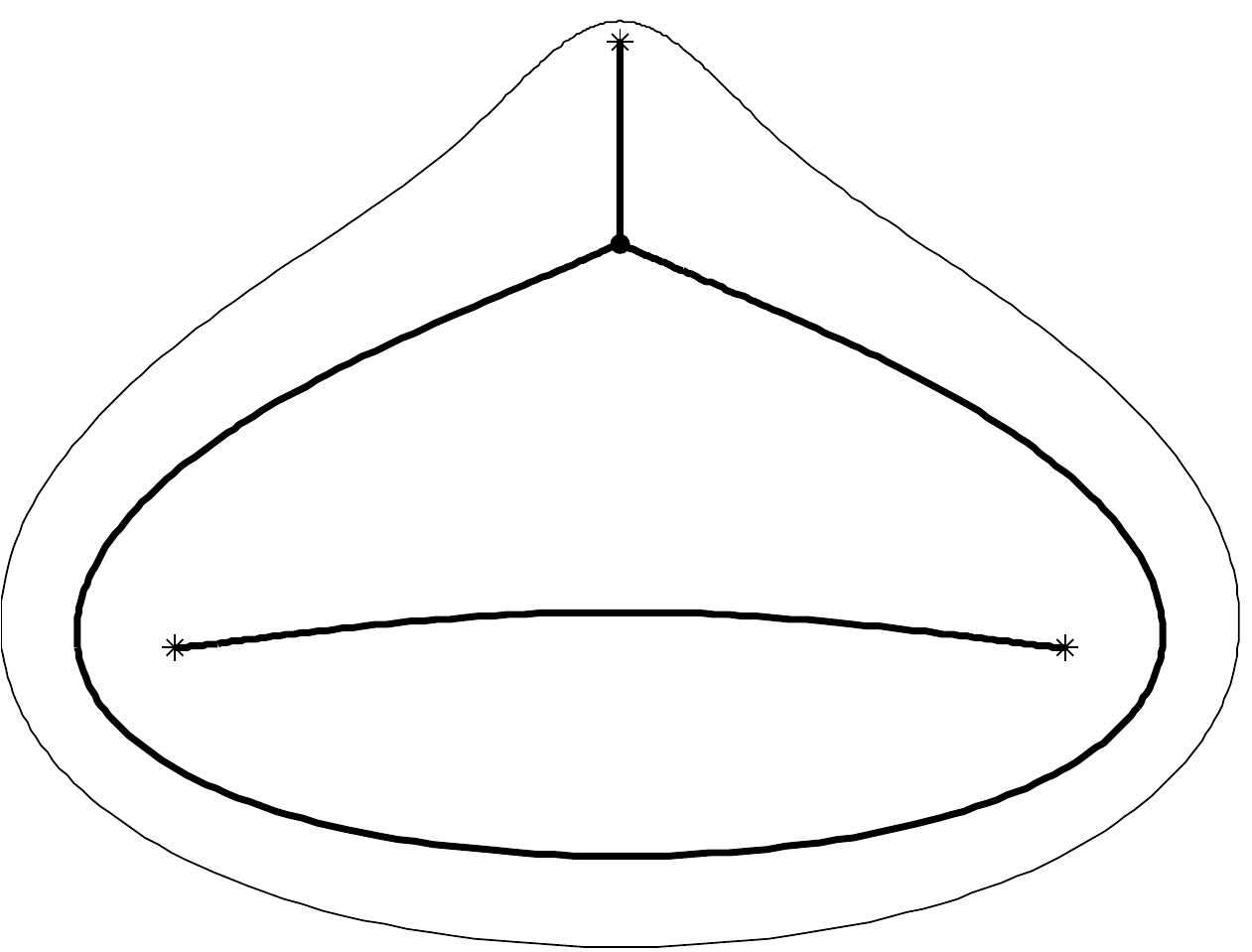}%
\end{overpic}}
 \\ 
\mbox{\begin{overpic}[scale=0.3]%
{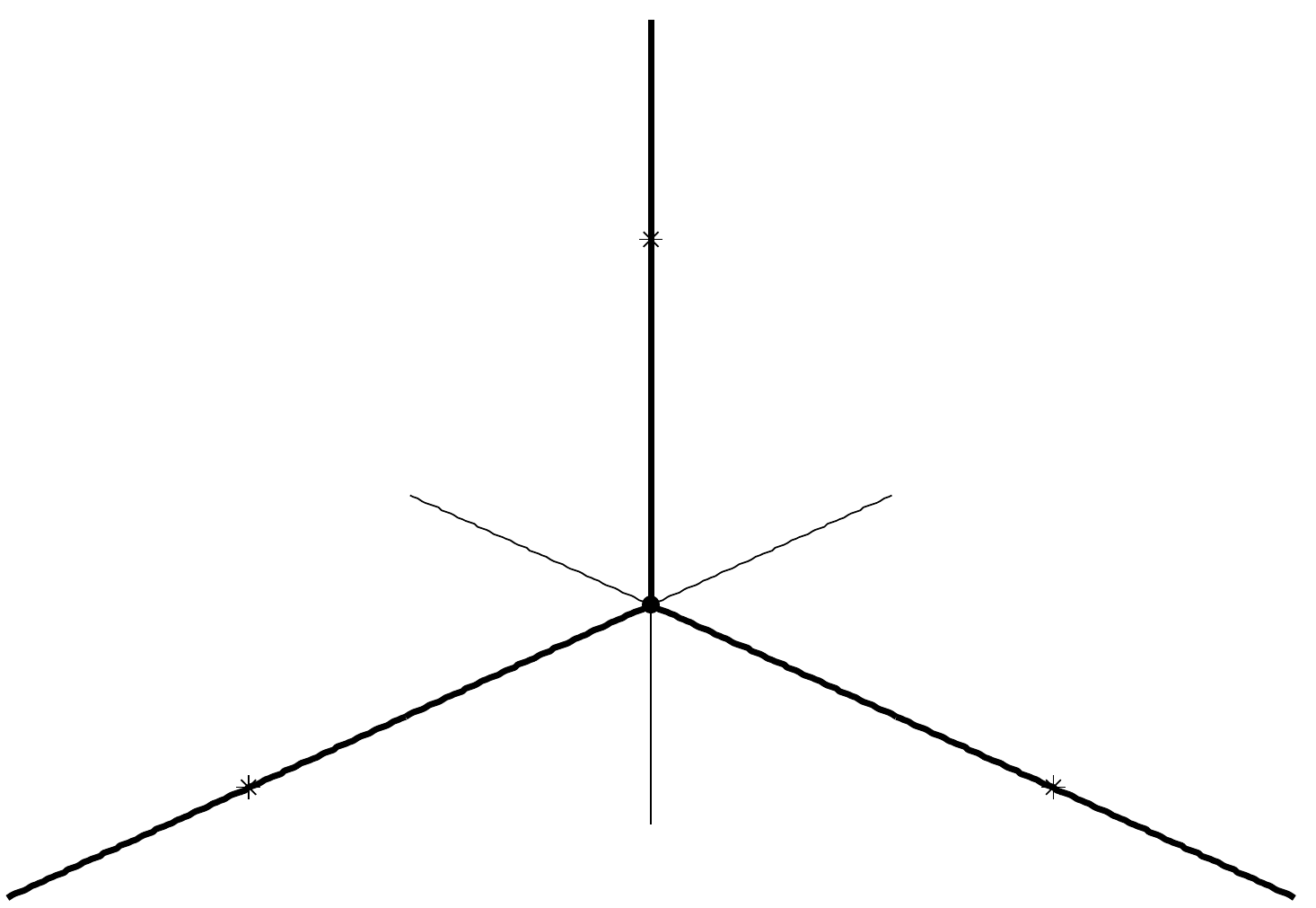}%
\put(58,45){\emph{iv)} $\mu_v(\gamma(v))=m_k$}
      \put(40,49){\small $a_k $}
      \put(23,20){\small $v=v^* $}
\end{overpic}} &
\mbox{\begin{overpic}[scale=0.35]%
{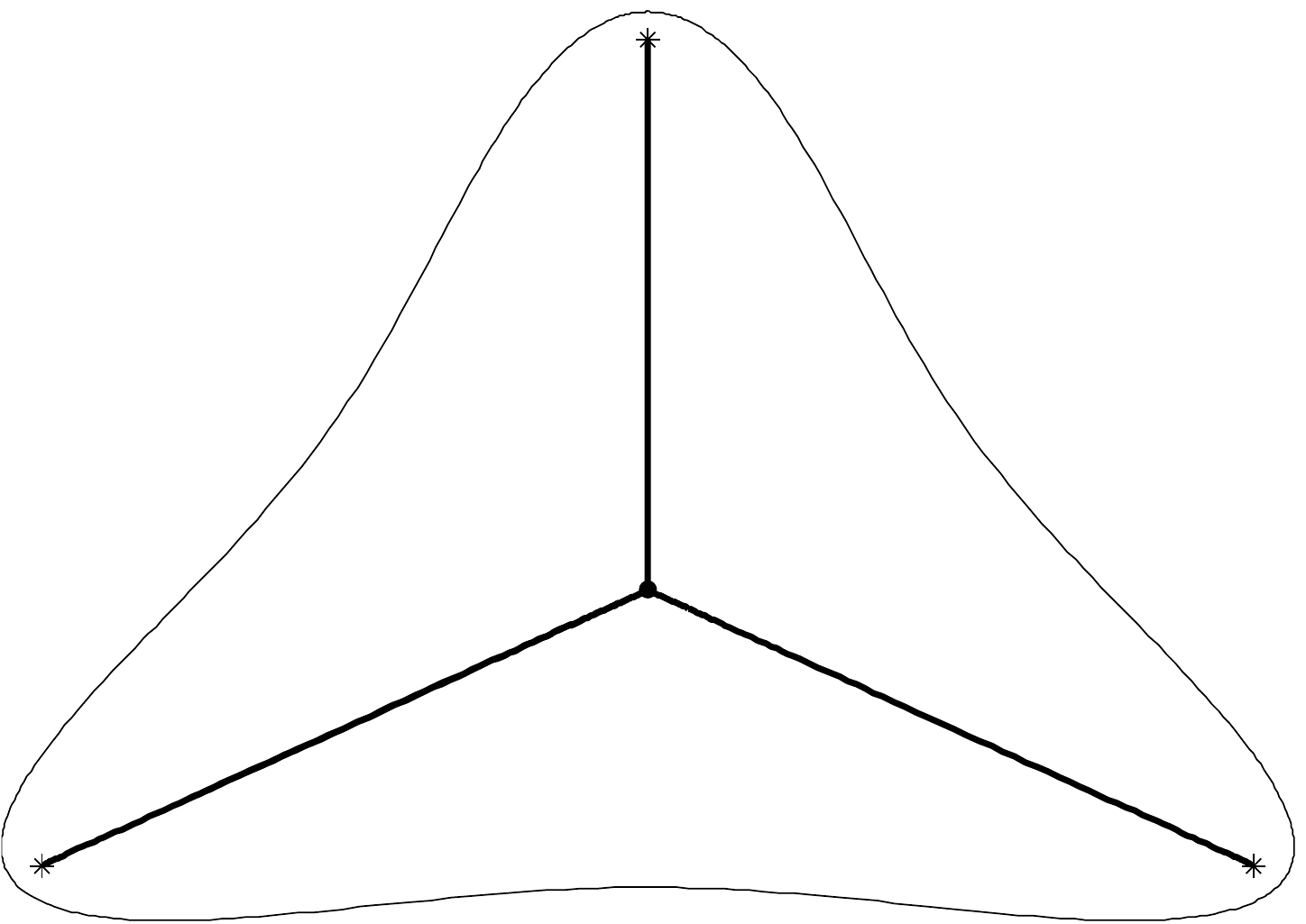}%
\end{overpic}} \\ 
\mbox{\begin{overpic}[scale=0.3]%
{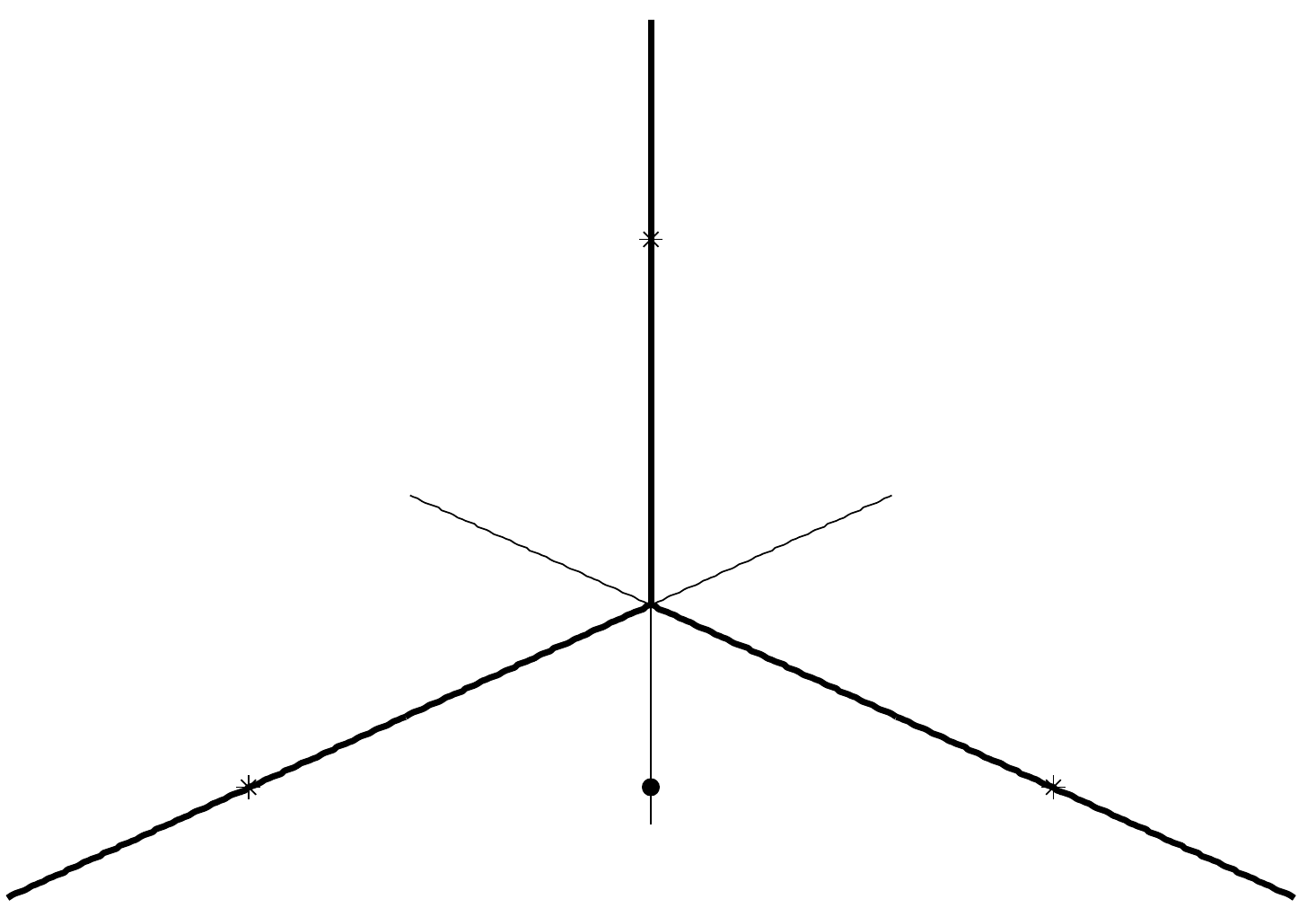}%
\put(58,45){\emph{v)} $\mu_v(\gamma'(v))<0$}
      \put(40,49){\small $a_k $}
      \put(43,8){\small $v $}
\end{overpic}} &
\mbox{\begin{overpic}[scale=0.35]%
{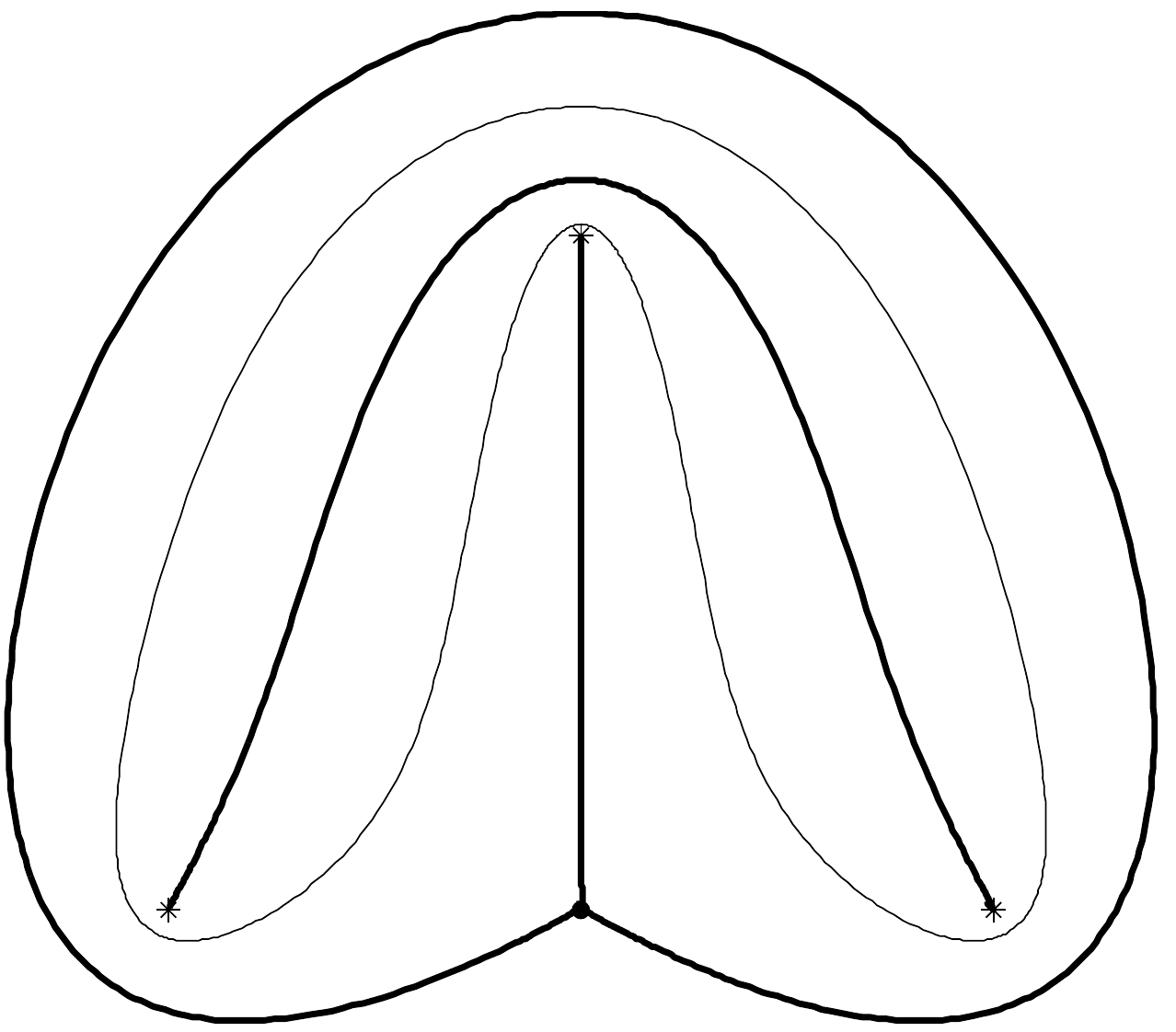}%
\end{overpic}}
\end{tabular}
\caption{Position of $v$ on $\ell_0 \cup \ell_1 \cup \ell_2$ (left) and the corresponding trajectories of the differential $\varpi_v$ in \eqref{qdp2}.}
\label{fig:cases}
\end{figure}

Furthermore, by Theorem \ref{thm:structureVPlus2}, the limit zero $v$ of $V(z)=z-v$ satisfies an equation of the form
$$
\frac{1}{\pi i}\, \int_{a_{k}}^{v} \sqrt{\frac{V(t)}{A(t)}}\, dt =  \beta \in [0, M_{k}],
$$
which is consistent with our construction in Subsection \ref{subs:p=2Cheb}. We can go back now to  the set of equations \eqref{VVequationsGeneral} and add the following result:
 \begin{theorem}
For any $\varepsilon>0$, equations
\begin{equation}\label{quasiVV}
w_{k}(\widetilde v_{j,k})= \frac{ j }{\lambda_{n}}+ \frac{ \rho_{k}/2  }{\lambda_{n}} , \quad j=0, \dots, M_{k} = [m_{k}n(1-\varepsilon)],
\end{equation}
with $k=0, 1, 2$, uniquely define $\widetilde n=(1-\varepsilon) n\pm 1$ points $\widetilde v_{j,k}$. There exists a constant $C=C(\varepsilon)>0$ such that any Van Vleck zero lies in a $C/n^2$ distance from a point  $\widetilde v_{j,k}$.
\end{theorem}
In this way, each zero of $V_{n}$ ``belongs'' to a $\mathcal O(n^{-2})$ neighborhood of one of the points defined by the equations above. We need to show additionally that (again, up to a small neighborhood of the Chebotarev's center $v^*$) all neighborhoods of this form contain at least one Van Vleck zero (establishing in this form a bijection between points $\widetilde v_{j,k}$ and the set $\mathcal V_{n}$). However, the techniques developed in this paper do not allow to complete this proof, so we formulate it as an open question:
\begin{conjecture}
There exists a constant $C=C(\varepsilon)>0$ such that to each point $\widetilde v_{j,k}$ defined by equations \eqref{quasiVV} it corresponds at least one zero $v$ of a Van Vleck polynomial such that $|v-\widetilde v_{j,k}|\leq C/n^2$.\end{conjecture}

\section*{Acknowledgments}

The authors gratefully acknowledge the help of Dar\'{\i}o Ramos-L\'opez with the numerical experiments yielding Figure \ref{fig:experiments}.

\def\cprime{$'$}
\providecommand{\bysame}{\leavevmode\hbox to3em{\hrulefill}\thinspace}
\providecommand{\MR}{\relax\ifhmode\unskip\space\fi MR }
\providecommand{\MRhref}[2]{%
  \href{http://www.ams.org/mathscinet-getitem?mr=#1}{#2}
}
\providecommand{\href}[2]{#2}

\end{document}